\documentclass[11pt,draft]{article}
\usepackage{amsmath,amscd}
\usepackage{amssymb,latexsym,amsthm}
\usepackage{color}
\usepackage[spanish,english]{babel}
\usepackage{amssymb}
\usepackage{color}
\usepackage{amsmath,amsthm,amscd}
\usepackage[latin1]{inputenc}
\usepackage{lscape}
\usepackage{amsfonts}

\numberwithin{equation}{section}
\newtheorem{theorem}{Theorem}[section]

\newtheorem{proposition}[theorem]{Proposition}

\newtheorem{corollary}[theorem]{Corollary}

\theoremstyle{definition}

\newtheorem{definition}[theorem]{Definition}
\newtheorem{examples}[theorem]{Examples}

\newtheorem{remark}[theorem]{Remark}

\newcommand{\cA}{\mbox{${\cal A}$}}

\newcommand{\cU}{\mbox{${\cal U}$}}

\title{\textbf{Armendariz Modules over Skew PBW Extensions}}
\author{Armando Reyes\footnote{Departamento de  Matem\'aticas. e-mail: mareyesv@unal.edu.co} \\ Universidad Nacional de Colombia, sede Bogot\'a}
\date{}
\begin{document}
\maketitle
\begin{abstract}
\noindent The aim of this paper is to develop the theory of skew Armendariz and quasi-Armendariz modules over skew PBW extensions. We generalize the results of several works in the literature concerning Ore extensions to another non-commutative rings which can not be expressed as iterated Ore extensions. As a consequence of our treatment, we extend and unify different results about the Armendariz, Baer, p.p., and p.q.-Baer properties for Ore extensions and skew PBW extensions.

\bigskip


\bigskip

\bigskip

\end{abstract}
\section{Introduction}\label{section}
In \cite{Kaplansky1968}, Kaplansky defined a ring $B$ as a {\em Baer} ({\em quasi-Baer}, which was defined by Clark \cite{Clark1967}) ring, if the right annihilator of every nonempty subset (ideal) of $B$ is generated by an idempotent (the objective of these rings is to abstract various properties of von Neumann algebras and complete $\ast$-regular rings; Clark used the quasi-Baer concept to characterize when a finite-dimensional algebra with unity over an algebraically closed field is isomorphic to a twisted matrix units semigroup algebra). Ano\-ther generalization of Baer rings are the p.p.-rings. A ring $B$ is called {\em right} ({\em left}) {\em principally projective} ({\em p.p.} for short), if the right ({\em left}) annihilator of each element of $B$ is ge\-ne\-ra\-ted by an idempotent (or equivalently, rings in which each principal right ({\em left}) ideal is projective).  Birkenmeier et al. \cite{Birkenmeieretal2001a} defined a ring to be called a {\em right} ({\em left}) {\em principally quasi-Baer} (or simply {\em right} ({\em left}) {\em p.q.-Baer}) ring, if the right annihilator of each principal right (left) ideal of $B$ is generated by an idempotent. Note that in a reduced ring $B$, $B$ is Baer (p.p.) if and only if $B$ is quasi-Baer (p.q.-Baer).\\

Commutative and noncommutative  Baer, quasi-Baer, p.p. and p.q.-Baer  rings have been investigated in the literature. For instance, in  \cite{Armendariz1974},  Armendariz established the following proposition: if $B$ is a reduced ring, then $B[x]$ is a Baer ring if and only if $B$ is a Baer ring (\cite{Armendariz1974}, Theorem B). In fact, Armendariz showed an example to illustrate that the condition to be reduced is not superfluous. Birkenmeier et. al., in \cite{Birkenmeieretal2001a}  showed that the quasi-Baer condition is preserved by many polynomial extensions, and in \cite{Birkenmeieretal2000}, they proved that a ring $B$ is right p.q.-Baer if and only if $B[x]$ is right p.q.-Baer. In the context of Ore extensions (defined by Ore in \cite{Ore1933}) given by $B[x;\sigma,\delta]$ with $\sigma$ injective (also known as Ore extensions of injective type), we found  several works in the literature, see  \cite{Clark1967}, \cite{Birkenmeieretal2000}, \cite{HongKimKwak2000}, \cite{Birkenmeieretal2001a},  \cite{BirkenmeierKimPark2001},  \cite{Hirano2002},  \cite{HashemiMoussavi2005}, and others (in \cite{Reyes2015}, \cite{ReyesSuarezUMA2017}, or \cite{ReyesSuarezClifford2017} a detailed list of references is presented). Some of these works consider the case $\delta=0$ and $\sigma$ an automorphism, or the case where $\sigma$ is the identity. It is important to say that the Baerness and quasi-Baerness of a ring $B$ and an Ore extension $B[x;\sigma,\delta]$ of $B$ does not depend on each other. More exactly, there are examples which show that there exists a Baer ring $B$ but the Ore extension $B[x;\sigma,\delta]$ is not right p.q.-Baer; similarly, there exist Ore extensions $B[x;\sigma,\delta]$ which are quasi-Baer, but $B$ is not quasi-Baer (see \cite{HongKimKwak2000}, Examples 8, 9 and 10 for more details). \\

With respect to the context of modules, Lee and Zhou in \cite{LeeZhou2004} introduced the notions of Baer, quasi-Baer and p.p.-modules in the following way: for a ring $B$ and a right $B$-module $M_B$, (i) $M_B$ is called {\em Baer} ({\em quasi-Baer}) if, for any subset (submodule) $X$ of $M$, ${\rm ann}_B(X) = eB$, where $e^{2}=e\in B$; (ii) $M_B$ is called {\em principally projective} ({\em p.p.} for short) {\em module} ({\em principally quasi-Baer} {\em module}) if, for any element $m\in M$, ${\rm ann}_B(m) = eB$ (${\rm ann}_B(mB) = eB$), where $e^{2}=e\in B$. It is important to remark that all these notions coincide with the ring definitions above, considering a ring $B$ as a right module over itself. In other words, a ring $B$ is Baer (quasi-Baer or p.p.) if and only if $B_B$ is a Baer (quasi-Baer or p.p.) module. In fact, if $B$ is a Baer (quasi-Baer or p.p.) ring, then for any right ideal $I$ of $B$, $I_B$ is a Baer (quasi-Baer or p.p.) module. Note that $B$ is a right p.q.-Baer ring if and only if $B_B$ is a p.q.-Baer module, and every submodule of a p.q.-Baer module is p.q.-Baer, and every Baer module is quasi-Baer.\\

Since the notion of reduced ring (a ring $B$ is called {\em reduced} if it has no nonzero nilpotent elements; note that every reduced ring is abelian, i.e., every idempotent is central) is very important for characterizing the properties of being p.p. and p.q.-Baer (in \cite{HongKimKwak2000}, Lemma 1 it was proved that for a reduced ring $B$,  $B$ is a right p.p.-ring $\Leftrightarrow B$ is a p.p.-ring $\Leftrightarrow B$ is a right p.q.-Baer ring $\Leftrightarrow B$ is a p.q.-Baer ring), it is of interest to know its corresponding notion for the context of modules: $M_B$ is called {\em reduced} (Lee and Zhou \cite{LeeZhou2004}), if for any elements $m\in M,\ a\in B$, $ma=0$ implies $mB\cap Ma = 0$. Precisely, Lee and Zhou generalized several results of reduced rings to reduced modules.\\

The notion of Armendariz ring, which is the primary object of study in this paper, it has also been investigated. Let us recall briefly. In commutative algebra, a ring $B$ is called {\em Armendariz} (the term was introduced by Rege and Chhawchharia in \cite{RegeChhawchharia1997}) if whenever polynomials $f(x)=a_0+a_1x+\dotsb + a_nx^n$, $g(x)=b_0+b_1x+\dotsb + b_mx^m\in B[x]$ satisfy $f(x)g(x)=0$, then $a_ib_j=0$, for every $i,j$. The interest of this notion lies in its natural and its useful role in understanding the relation between the annihilators of the ring $B$ and the annihilators of the polynomial ring $B[x]$. In  \cite{Armendariz1974}, Lemma 1, Armendariz showed that a  reduced ring always satisfies this condition. Now, in the context of Ore extensions, Armendariz property has also been studied. For instance, Hirano in \cite{Hirano2002} defined a ring $B$ to be {\em quasi-Armendariz} if whenever two polynomials $f(x)=\sum_{i=0}^{m} a_ix^{i},\ g(x) = \sum_{j=0}^{t} b_jx^{j}\in B[x]$ satisfy $f(x)B[x]g(x)=0$, then $a_iRb_j$, for every $i, j$.
In \cite{HongKimKwak2003}, Hong et. al., extended the Armendariz property of rings to skew polynomial rings $B[x;\alpha]$ with zero derivation. For an endomorphism $\alpha$ of a ring $B$, $B$ is called an $\alpha$-skew Armendariz ring, if for polynomials $f(x)=a_0 + a_1x + \dotsb + a_nx^{n}$ and $g(x) = b_0 + b_1x + \dotsb + b_mx^{m}$ in $B[x;\alpha]$, $f(x)g(x)=0$ implies $a_i\alpha^{i}(b_j)=0$, for every $0\le i\le n$, and $0\le j\le m$. A more general treatment for the notion of Armendariz for Ore extensions with $\delta$ not necessarily zero, it was established by Nasr-Isfahani and Moussavi \cite{NasrMoussavi2008} using the notion of skew-Armendariz ring. It is important to say that the relations between Armendariz rings and Baer (quasi-Baer) rings have been also investigated in different papers, see for example \cite{Armendariz1974}, \cite{RegeChhawchharia1997}, \cite{AndersonCamillo1998}, \cite{Birkenmeieretal2000},   \cite{HongKimKwak2000},  \cite{Birkenmeieretal2001a},   \cite{Hirano2002},   \cite{HongKimKwak2003},  \cite{Matczuk2004}, and others (see \cite{Reyes2015}, \cite{ReyesSuarezUMA2017}, or \cite{ReyesSuarezClifford2017} for a detailed list of references). \\

The notion of Armendariz for modules over Ore extensions also have been formulated. In \cite{ZhangChen2008}, Zhang and Chen introduced the notion of $\alpha$-skew Armendariz modules over Ore extensions with zero derivation ($\delta=0$) in the following way: an $B$-module $M$ is called $\alpha$-{\rm skew Armendariz}, if for polynomials $m(x)=m_0 + m_1x + \dotsb + m_kx^{k}\in M[X]$ and $f(x) = b_0+b_1x +\dotsb + b_nx^{n}\in B[x;\alpha]$, $m(x)f(x)=0$ implies $m_i\alpha^{i}(b_j)=0$, for every $0\le i\le k$ and $0\le j\le n$  (Baser in \cite{Baser2006} studied the relations between the set of annihilators in $M_B$ and the set of annihilators in $M[X]$). A module $M_B$ is called $\alpha$-{\em Armendariz}, if $M_B$ is $\alpha$-compatible and $\alpha$-skew-Armendariz (\cite{LeeZhou2004}). These authors also proved that $B$ is an $\alpha$-skew Armendariz ring if and only if every flat right $B$-module is $\alpha$-skew Armendariz, and a module $M_B$ is $\alpha$-reduced, if $M_B$ is $\alpha$-compatible and reduced. A more general treatment about the notion of Armendariz module over Ore extensions with $\delta$ not necessarily zero it was presented by Alhevaz and Moussavi in \cite{AlhevazMoussavi2012}. There, they study the relationship between an $B$-module $M_B$ and the general polynomial module $M[X]$ over the Ore extension $B[x;\alpha, \delta]$, and introduce the notions of skew-Armendariz modules and skew quasi-Armendariz modules which are generalizations of $\alpha$-skew Armendariz modules \cite{ZhangChen2008} and $\alpha$-reduced modules \cite{LeeZhou2004}. In fact, they also established several connections of the Baer, quasi-Baer and the p.p.-properties with the notion of skew Armendariz and skew quasi-Armendariz module. In this way, \cite{AlhevazMoussavi2012} extends and unifies several known results related to Armendariz rings and modules, such as \cite{HongKimKwak2000}, \cite{HongKimKwak2003}, \cite{NasrMoussavi2008}, \cite{ZhangChen2008}, and others, to general polynomial modules over Ore extensions.\\

With the aim of generalizing the results established about Armendariz and Baer properties in the mentioned papers above, in this article we are interested in a class of non-commutative rings of polynomial type more general than iterated Ore extensions (of injective type), the {\em skew Poincar\'e-Birkhoff-Witt extensions} (also known as $\sigma$-{\em PBW} {\em extensions}), where PBW denotes Poincar\'e-Birkhoff-Witt, introduced in \cite{LezamaGallego2011} (see Examples \ref{mentioned} for a list of non-commutative rings which are $\sigma$-PBW extensions but not iterated Ore extensions). Actually, skew PBW extensions are more general than several families of non-commutative rings, such as universal enveloping algebras of finite dimensional Lie algebras, PBW extensions introduced by Bell and Goodearl in \cite{BellGoodearl1988}, almost normalizing extensions defined by McConnell and Robson in \cite{McConnellRobson2001}, solvable polynomial rings introduced by Kandri-Rody and Weispfenning in \cite{KandryWeispfenninig1990}, and generalized by Kredel in \cite{Kredel1992},  diffusion algebras studied by Isaev, Pyatov, and Rittenberg  in \cite{IsaevPyatovRittenberg2001}, and other kind of non-commutative algebras of polynomial type. The importance of skew PBW extensions is that the coefficients do not necessarily commute with the variables, and these  coefficients are not necessarily elements of fields (see Definition \ref {gpbwextension} below). In fact, the $\sigma$-PBW extensions contain well-known groups of algebras such as some types of $G$-algebras studied by Levandovskyy \cite{Levandovskyy2005} and some PBW algebras defined by Bueso et. al., in \cite{BuesoGomezTorrecillas2003} (both $G$-algebras and PBW algebras take coefficients in fields and assume that coefficientes commute with variables),  Auslander-Gorenstein rings, some Calabi-Yau and skew Calabi-Yau algebras, some Artin-Schelter regular algebras, some Koszul algebras, quantum polynomials, some quantum universal enveloping algebras, and others (see \cite{Reyes2013PhD},
\cite{LezamaReyes2014}, and \cite{SuarezLezamaReyes2015}  for a detailed list of examples). For more details about the relation between $\sigma$-PBW extensions and another algebras with PBW bases, see \cite{Reyes2013PhD} or \cite{LezamaReyes2014}.\\

Since Ore extensions of injective type are particular examples of $\sigma$-PBW extensions, and ha\-ving in mind that several ring, module and homological properties of have been studied by the author and others for skew PBW extensions (see \cite{LezamaGallego2011},  \cite{Reyes2013PhD}, \cite{Reyes2013}, \cite{Reyes2014}, \cite{Reyes2014b},   \cite{Artamonov2015},  \cite{LezamaAcostaReyes2015},  \cite{Reyes2015}, \cite{ArtamonovLezamaFajardo2016}, \cite{ReyesSuarezUMA2017},  \cite{ReyesSuarezClifford2017},  \cite{ReyesSuarezUPTC2016}, etc), we consider relevant to investigate the properties of Baer, quasi-Baer, p.p., p.q.-Baer, and Armendariz in the context of modules over these extensions (in \cite{Reyes2015}, \cite{ReyesSuarez2016},  \cite{ReyesSuarez2016b},  \cite{ReyesSuarezUMA2017}, and \cite{ReyesSuarezClifford2017},    these properties were investigated for $\sigma$-PBW extensions, for example, with the purpose of computing its Goldie dimension \cite{Reyes2014}) with the aim of establishing and generalizing several results in the literature for Ore extensions of injective type and $\sigma$-PBW extensions. In this way, our results generalizes several works concerning Ore extensions and $\sigma$-PBW extensions, such as \cite{HongKimKwak2000},  \cite{Hirano2002}, \cite{HongKimKwak2003},  \cite{LeeZhou2004}, \cite{MoussaviHashemi2005}, \cite{NasrMoussavi2008}, \cite{ZhangChen2008}, \cite{AlhevazMoussavi2012}, \cite{Reyes2015}, \cite{ReyesSuarez2016b}, \cite{ReyesSuarezClifford2017}, and \cite{ReyesSuarezUMA2017}. We can say that the importance of our results is precisely to establish all these properties for those non-commutative rings which can not be expressed as Ore extensions.\\

The paper is organized as follows: In Section \ref{definitionexamplesspbw} we establish some useful results about $\sigma$-PBW extensions for the rest of the paper. Next, in Section \ref{armendarizmodules} we introduce the notion of skew Armendariz and skew quasi-Armendariz modules based on \cite{AlhevazMoussavi2012}. First, Section \ref{skewarmendarizmodules} contains the definition of skew-Armendariz module for $\sigma$-PBW extensions, and in Section \ref{SkewquasiArmendarizmodules} we introduce the notion of skew quasi-Armendariz module. The more important results of this paper are presented in this section following the ideas established by Alhevaz and Moussavi in \cite{AlhevazMoussavi2012} for the case of Ore extensions. It is a remarkable fact that the tools employed in that paper are very useful for the study of Armendariz modules over rings which can not be expressed as Ore extensions. In this way, the techniques used here are fairly standard and follow the same path as other text on the subject. The results presented are new for skew PBW extensions and all they generalize  others existing in the literature.\\

Throughout the paper, the word ring means a ring not necessarily commutative with unity, and all modules are right modules and ${\rm ann}_B(X):=\{r\in B\mid Xr=0\}$, where $X\subseteq B$, for any ring $B$.
\section{Skew PBW extensions}\label{definitionexamplesspbw}
In this section we establish some useful results about skew PBW extensions for the rest of the paper.
\begin{definition}[\cite{LezamaGallego2011}, Definition 1]\label{gpbwextension}
Let $R$ and $A$ be rings. We say that $A$ is  a {\em $\sigma$-PBW extension} (also known as {\em skew PBW extension}) {\em of}  $R$, which is denoted by $A:=\sigma(R)\langle
x_1,\dots,x_n\rangle$, if the following conditions hold:
\begin{enumerate}
\item[\rm (i)]$R\subseteq A$;
\item[\rm (ii)]there exist elements $x_1,\dots ,x_n\in A$ such that $A$ is a left free $R$-module, with basis ${\rm Mon}(A):= \{x^{\alpha}=x_1^{\alpha_1}\cdots
x_n^{\alpha_n}\mid \alpha=(\alpha_1,\dots ,\alpha_n)\in
\mathbb{N}^n\}$,  and $x_1^{0}\dotsb x_n^{0}:=1\in {\rm Mon}(A)$.

\item[\rm (iii)]For each $1\leq i\leq n$ and any $r\in R\ \backslash\ \{0\}$, there exists an element $c_{i,r}\in R\ \backslash\ \{0\}$ such that $x_ir-c_{i,r}x_i\in R$.
\item[\rm (iv)]For any elements $1\leq i,j\leq n$, there exists $c_{i,j}\in R\ \backslash\ \{0\}$ such that $x_jx_i-c_{i,j}x_ix_j\in R+Rx_1+\cdots +Rx_n$.
\end{enumerate}
\end{definition}
\begin{proposition}[\cite{LezamaGallego2011}, Proposition
3]\label{sigmadefinition}
Let $A$ be a $\sigma$-PBW  extension of $R$. For each $1\leq i\leq
n$, there exist an injective endomorphism $\sigma_i:R\rightarrow
R$ and an $\sigma_i$-derivation $\delta_i:R\rightarrow R$ such that $x_ir=\sigma_i(r)x_i+\delta_i(r)$, for  each $r\in R$. We write  $\Sigma:=\{\sigma_1,\dotsc, \sigma_n\}$, and $\Delta:=\{\delta_1,\dotsc, \delta_n\}$.
\end{proposition}
\begin{definition}[\cite{LezamaGallego2011}, Definition 4]\label{sigmapbwderivationtype}
Let $A$ be a $\sigma$-PBW extension of $R$.
\begin{enumerate}
\item[\rm (a)] $A$ is called \textit{quasi-commutative} if the conditions
{\rm(}iii{\rm)} and {\rm(}iv{\rm)} in Definition
\ref{gpbwextension} are replaced by the following: (iii') for each $1\leq i\leq n$ and all $r\in R\ \backslash\ \{0\}$, there exists $c_{i,r}\in R\ \backslash\ \{0\}$ such that $x_ir=c_{i,r}x_i$; (iv') for any $1\leq i,j\leq n$, there exists $c_{i,j}\in R\ \backslash\ \{0\}$ such that $x_jx_i=c_{i,j}x_ix_j$.
\item[\rm (b)] $A$ is called \textit{bijective}, if $\sigma_i$ is bijective for each $1\leq i\leq n$, and $c_{i,j}$ is invertible, for any $1\leq
i<j\leq n$.
\end{enumerate}
\end{definition}
\begin{examples}\label{mentioned}
If $R[x_1;\sigma_1,\delta_1]\dotsb [x_n;\sigma_n,\delta_n]$ is an iterated Ore extension where
\begin{itemize}
\item $\sigma_i$ is injective, for $1\le i\le n$;
\item $\sigma_i(r)$, $\delta_i(r)\in R$, for every $r\in R$ and $1\le i\le n$;
\item $\sigma_j(x_i)=cx_i+d$, for $i < j$, and $c, d\in R$, where $c$ has a left inverse;
\item $\delta_j(x_i)\in R + Rx_1 + \dotsb + Rx_n$, for $i < j$,
\end{itemize}
then $R[x_1;\sigma_1,\delta_1]\dotsb [x_n;\sigma_n, \delta_n] \cong \sigma(R)\langle x_1,\dotsc, x_n\rangle$ (\cite{LezamaReyes2014}, p. 1212). Note that $\sigma$-PBW extensions of endomorphism type are more general than iterated Ore extensions $R[x_1;\sigma_1]\dotsb [x_n;\sigma_n]$, and in general, $\sigma$-PBW extensions are more general than Ore extensions of injective type.

Next, we present some non-commutative rings which are $\sigma$-PBW extensions but they can not be expressed as iterated Ore extensions (see \cite{LezamaReyes2014} for the reference of every example).
\begin{enumerate}
\item[\rm (a)] Let $k$ be a
commutative ring and $\mathfrak{g}$ a finite dimensional Lie
algebra over $k$ with basis $\{x_1,\dots ,x_n\}$. The \textit{universal
enveloping algebra} of
$\mathfrak{g}$, denoted $\cU(\mathfrak{g})$, is a skew PBW extension
of $k$, since $x_ir-rx_i=0$, $x_ix_j-x_jx_i=[x_i,x_j]\in
\mathfrak{g}=k+kx_1+\cdots+kx_n$, $r\in k$, for $1\leq i,j\leq n$. In particular, the \textit{universal enveloping algebra} \textit{of a Kac-Moody Lie algebra} is a skew PBW extension
of a polynomial ring.
\item [\rm (b)] The \textit{universal enveloping ring} $\cU(V,R
,\Bbbk)$, where  $R$ is a
$\Bbbk$-algebra, and $V$ is a $\Bbbk$-vector space which is also a
Lie ring containing $R$ and $\Bbbk$ as Lie ideals with suitable
relations. The enveloping ring $\cU(V,R,\Bbbk)$ is a finite skew
 PBW extension of $R$ if ${\rm dim}_\Bbbk\ (V/R)$ is finite.
\item [\rm (c)] Let $k$, $\mathfrak{g}$, $\{x_1,\dots ,x_n\}$ and $\cU(\mathfrak{g})$ be as in the
previous example; let $R$ be a $k$-algebra containing $k$. The
\textit{tensor product} $A:=R\ \otimes_k\ \cU(\mathfrak{g})$ is a skew PBW extension of $R$, and it is a particular case
of \textit{crossed product} $R*\cU(\mathfrak{g})$ of $R$ by $\cU(\mathfrak{g})$, which is a skew PBW  extension of $R$.
\item [\rm (d)] The \textit{twisted or smash product differential operator ring} $R\ \# _{\sigma}\  \cU(\mathfrak{g})$, where $\mathfrak{g}$ is a finite-dimensional Lie algebra acting on
$R$ by derivations, and $\sigma$ is Lie 2-cocycle with values in
$R$.
\item [\rm (e)] Diffusion algebras arise in physics as a possible way to
understand a large class of $1$-dimensional stochastic process \cite{IsaevPyatovRittenberg2001}. A \textit{diffusion algebra} $\cA$ with
parameters $a_{ij}\in \mathbb{C}\ \backslash\ \{0\}, 1\le i, j\le n$, is an
algebra over $\mathbb{C}$ generated by variables $x_1,\dotsc,x_n$
subject to relations $a_{ij}x_ix_j-b_{ij}x_jx_i=r_jx_i-r_ix_j$, whenever $i<j$, $b_{ij}, r_i\in \mathbb{C}$ for all $i<j$. $\cA$ admits a $PBW$-basis of standard monomials $x_1^{i_1}\dotsb x_n^{i_n}$, that is, $\cA$
is a diffusion algebra if these standard monomials are a $\mathbb{C}$-vector space basis for $\cA$. From Definition \ref{gpbwextension}, (iii)
and (iv), it is clear that the family of skew PBW extensions are more general than diffusion algebras.  We will denote
$q_{ij}:=\frac{b_{ij}}{a_{ij}}$. The parameter $q_{ij}$ can be a
root of unity if and only if is equal to 1. It is therefore reasonable to assume that these parameters not to be a root of unity other
than 1. If all coefficients $q_{ij}$ are
nonzero, then the corresponding diffusion algebra have a PBW
basis of standard monomials $x_1^{i_1}\dotsb x_n^{i_n}$, and hence
these algebras are skew PBW extensions. More precisely,
$\cA\cong \sigma(\boldsymbol{\mathbb{C}})\langle
x_1,\dotsc,x_n\rangle$.
\end{enumerate}
It is important to say that $\sigma$-PBW extensions contains various well-known groups of algebras such as PBW extensions \cite{BellGoodearl1988}, the almost normalizing extensions  \cite{McConnellRobson2001}, solvable polynomial rings \cite{KandryWeispfenninig1990}, and  \cite{Kredel1992}, diffusion algebras \cite{IsaevPyatovRittenberg2001}, some types of Auslander-Gorenstein rings, some skew Calabi-Yau algebras, some Artin-Schelter regular algebras, some Koszul algebras, quantum polynomials, some quantum universal enveloping algebras,  etc. In comparison with $G$-algebras \cite{Levandovskyy2005} or PBW algebras \cite{BuesoGomezTorrecillas2003}, $\sigma$-PBW extensions do not assume that the ring of coefficients is a field neither that the coefficients commute with the variables, so that skew PBW extensions are not included in these algebras. Indeed, the $G$-algebras with  $d_{i,j}$ linear (recall that for these algebras $x_jx_i = c_{i,j}x_ix_j+ d_{i,j},\ 1\le i < j \le n$), are particular examples of $\sigma$-PBW extensions.  A detailed list of examples of skew PBW extensions and its relations with another algebras with PBW bases is presented in \cite{Reyes2013PhD} and \cite{LezamaReyes2014}.
\end{examples}
\begin{definition}[\cite{LezamaGallego2011}, Definition 6]\label{definitioncoefficients}
Let $A$ be a $\sigma$-PBW extension of $R$. Then:
\begin{enumerate}
\item[\rm (i)]for $\alpha=(\alpha_1,\dots,\alpha_n)\in \mathbb{N}^n$,
$\sigma^{\alpha}:=\sigma_1^{\alpha_1}\cdots \sigma_n^{\alpha_n}$,
$|\alpha|:=\alpha_1+\cdots+\alpha_n$. If
$\beta=(\beta_1,\dots,\beta_n)\in \mathbb{N}^n$, then
$\alpha+\beta:=(\alpha_1+\beta_1,\dots,\alpha_n+\beta_n)$.
\item[\rm (ii)]For $X=x^{\alpha}\in {\rm Mon}(A)$,
$\exp(X):=\alpha$, $\deg(X):=|\alpha|$, and $X_0:=1$. The symbol $\succeq$ will denote a total order defined on ${\rm Mon}(A)$ (a total order on $\mathbb{N}^n$). For an
 element $x^{\alpha}\in {\rm Mon}(A)$, ${\rm exp}(x^{\alpha}):=\alpha\in \mathbb{N}^n$.  If
$x^{\alpha}\succeq x^{\beta}$ but $x^{\alpha}\neq x^{\beta}$, we
write $x^{\alpha}\succ x^{\beta}$. Every element $f\in A$ can be expressed uniquely as $f=a_0 + a_1X_1+\dotsb +a_mX_m$, with $a_i\in R$, and $X_m\succ \dotsb \succ X_1$ (eventually, we will use expressions as $f=a_0 + a_1Y_1+\dotsb +a_mY_m$, with $a_i\in R$, and $Y_m\succ \dotsb \succ Y_1$). With this notation, we define ${\rm
lm}(f):=X_m$, the \textit{leading monomial} of $f$; ${\rm
lc}(f):=a_m$, the \textit{leading coefficient} of $f$; ${\rm
lt}(f):=a_mX_m$, the \textit{leading term} of $f$; ${\rm exp}(f):={\rm exp}(X_m)$, the \textit{order} of $f$; and
 $E(f):=\{{\rm exp}(X_i)\mid 1\le i\le t\}$. Note that $\deg(f):={\rm max}\{\deg(X_i)\}_{i=1}^t$. Finally, if $f=0$, then
${\rm lm}(0):=0$, ${\rm lc}(0):=0$, ${\rm lt}(0):=0$. We also
consider $X\succ 0$ for any $X\in {\rm Mon}(A)$. For a detailed description of monomial orders in skew PBW  extensions, see \cite{LezamaGallego2011}, Section 3.
\end{enumerate}
\end{definition}
\begin{proposition}[\cite{LezamaGallego2011}, Theorem 7]\label{coefficientes}
If $A$ is a polynomial ring with coefficients in $R$ with respect to the set of indeterminates $\{x_1,\dots,x_n\}$, then $A$ is a skew PBW  extension of $R$ if and only if the following conditions hold:
\begin{enumerate}
\item[\rm (i)]for each $x^{\alpha}\in {\rm Mon}(A)$ and every $0\neq r\in R$, there exist unique elements $r_{\alpha}:=\sigma^{\alpha}(r)\in R\ \backslash\ \{0\}$, $p_{\alpha ,r}\in A$, such that $x^{\alpha}r=r_{\alpha}x^{\alpha}+p_{\alpha, r}$,  where $p_{\alpha ,r}=0$, or $\deg(p_{\alpha ,r})<|\alpha|$ if
$p_{\alpha , r}\neq 0$. If $r$ is left invertible,  so is $r_\alpha$.
\item[\rm (ii)]For each $x^{\alpha},x^{\beta}\in {\rm Mon}(A)$,  there exist unique elements $c_{\alpha,\beta}\in R$ and $p_{\alpha,\beta}\in A$ such that $x^{\alpha}x^{\beta}=c_{\alpha,\beta}x^{\alpha+\beta}+p_{\alpha,\beta}$, where $c_{\alpha,\beta}$ is left invertible, $p_{\alpha,\beta}=0$, or $\deg(p_{\alpha,\beta})<|\alpha+\beta|$ if
$p_{\alpha,\beta}\neq 0$.
\end{enumerate}
\end{proposition}
\begin{remark}\label{juradpr}
About Proposition \ref{coefficientes}, we have two observations:
\begin{enumerate}
\item [\rm (i)] (\cite{Reyes2015}, Proposition 2.9) If $\alpha:=(\alpha_1,\dotsc, \alpha_n)\in \mathbb{N}^{n}$ and $r\in R$, then
{\scriptsize{\begin{align*}
x^{\alpha}r = &\ x_1^{\alpha_1}x_2^{\alpha_2}\dotsb x_{n-1}^{\alpha_{n-1}}x_n^{\alpha_n}r = x_1^{\alpha_1}\dotsb x_{n-1}^{\alpha_{n-1}}\biggl(\sum_{j=1}^{\alpha_n}x_n^{\alpha_{n}-j}\delta_n(\sigma_n^{j-1}(r))x_n^{j-1}\biggr)\\
+ &\ x_1^{\alpha_1}\dotsb x_{n-2}^{\alpha_{n-2}}\biggl(\sum_{j=1}^{\alpha_{n-1}}x_{n-1}^{\alpha_{n-1}-j}\delta_{n-1}(\sigma_{n-1}^{j-1}(\sigma_n^{\alpha_n}(r)))x_{n-1}^{j-1}\biggr)x_n^{\alpha_n}\\
+ &\ x_1^{\alpha_1}\dotsb x_{n-3}^{\alpha_{n-3}}\biggl(\sum_{j=1}^{\alpha_{n-2}} x_{n-2}^{\alpha_{n-2}-j}\delta_{n-2}(\sigma_{n-2}^{j-1}(\sigma_{n-1}^{\alpha_{n-1}}(\sigma_n^{\alpha_n}(r))))x_{n-2}^{j-1}\biggr)x_{n-1}^{\alpha_{n-1}}x_n^{\alpha_n}\\
+ &\ \dotsb + x_1^{\alpha_1}\biggl( \sum_{j=1}^{\alpha_2}x_2^{\alpha_2-j}\delta_2(\sigma_2^{j-1}(\sigma_3^{\alpha_3}(\sigma_4^{\alpha_4}(\dotsb (\sigma_n^{\alpha_n}(r))))))x_2^{j-1}\biggr)x_3^{\alpha_3}x_4^{\alpha_4}\dotsb x_{n-1}^{\alpha_{n-1}}x_n^{\alpha_n} \\
+ &\ \sigma_1^{\alpha_1}(\sigma_2^{\alpha_2}(\dotsb (\sigma_n^{\alpha_n}(r))))x_1^{\alpha_1}\dotsb x_n^{\alpha_n}, \ \ \ \ \ \ \ \ \ \ \ \ \ \sigma_j^{0}:={\rm id}_R\ \ {\rm for}\ \ 1\le j\le n.
\end{align*}}}
\item [\rm (ii)] (\cite{Reyes2015}, Remark 2.10) Using (i), it follows that for the product $a_iX_ib_jY_j$, if $X_i:=x_1^{\alpha_{i1}}\dotsb x_n^{\alpha_{in}}$ and $Y_j:=x_1^{\beta_{j1}}\dotsb x_n^{\beta_{jn}}$, then
\begin{align*}
a_iX_ib_jY_j = &\ a_i\sigma^{\alpha_i}(b_j)x^{\alpha_i}x^{\beta_j} + a_ip_{\alpha_{i1}, \sigma_{i2}^{\alpha_{i2}}(\dotsb (\sigma_{in}^{\alpha_{in}}(b)))} x_2^{\alpha_{i2}}\dotsb x_n^{\alpha_{in}}x^{\beta_j} \\
+ &\ a_i x_1^{\alpha_{i1}}p_{\alpha_{i2}, \sigma_3^{\alpha_{i3}}(\dotsb (\sigma_{{in}}^{\alpha_{in}}(b)))} x_3^{\alpha_{i3}}\dotsb x_n^{\alpha_{in}}x^{\beta_j} \\
+ &\ a_i x_1^{\alpha_{i1}}x_2^{\alpha_{i2}}p_{\alpha_{i3}, \sigma_{i4}^{\alpha_{i4}} (\dotsb (\sigma_{in}^{\alpha_{in}}(b)))} x_4^{\alpha_{i4}}\dotsb x_n^{\alpha_{in}}x^{\beta_j}\\
+ &\ \dotsb + a_i x_1^{\alpha_{i1}}x_2^{\alpha_{i2}} \dotsb x_{i(n-2)}^{\alpha_{i(n-2)}}p_{\alpha_{i(n-1)}, \sigma_{in}^{\alpha_{in}}(b)}x_n^{\alpha_{in}}x^{\beta_j} \\
+ &\ a_i x_1^{\alpha_{i1}}\dotsb x_{i(n-1)}^{\alpha_{i(n-1)}}p_{\alpha_{in}, b}x^{\beta_j}.
\end{align*}
In this way, when we compute every summand of $a_iX_ib_jY_j$ we obtain products of the coefficient $a_i$ with several evaluations of $b_j$ in $\sigma$'s and $\delta$'s depending of the coordinates of $\alpha_i$.
\end{enumerate}
\end{remark}
\section{Armendariz modules over $\sigma$-PBW extensions}\label{armendarizmodules}
In this section we introduce the notions of skew Armendariz module and skew quasi-Armendariz module over $\sigma$-PBW extensions. We start defining the modules which we are going to study.\\

From Definition \ref{gpbwextension} we know that if $A$ is a $\sigma$-PBW extension of a ring $R$, then $A$ is a left free $R$-module. Now, Remark \ref{juradpr} (i) says us how to multiply elements of $R$ with elements of ${\rm Mon}(A)$, so that if we consider a right $R$-module $M_R$, we can consider the polynomial module $M\langle X\rangle_A$ over $A$. More precisely, as a set, the elements of $M\langle X\rangle_A$ are of the form $m_0 + m_1X_1 + \dotsb + m_tX_t$, $m_i\in M_R$ and $X_i\in {\rm Mon}(A)$,  for every $i$. If $\alpha:=(\alpha_1,\dotsc, \alpha_n)\in \mathbb{N}^{n}$ and $r\in R$, then the action of $A$ on these elements follow the rule established Remark \ref{juradpr} (ii). This fact is precisely because it suffices to define the action of monomials of $A$ on monomials in $M\langle X\rangle_A$. In other words, if $m_ix_1^{\alpha_{i1}}\dotsb x_n^{\alpha_{in}}$ and $b_jx_1^{\beta_{j1}}\dotsb x_n^{\beta_{jn}}$ are elements of $M\langle X\rangle$ and $A$, respectively, then we multiply these both elements following the rule
\begin{align}
m_ix_1^{\alpha_{i1}}\dotsb x_n^{\alpha_{in}}b_jx_1^{\beta_{j1}}\dotsb x_n^{\beta_{jn}} = &\ m_i\sigma^{\alpha_i}(b_j)x^{\alpha_i}x^{\beta_j} + m_ip_{\alpha_{i1}, \sigma_{i2}^{\alpha_{i2}}(\dotsb (\sigma_{in}^{\alpha_{in}}(b)))} x_2^{\alpha_{i2}}\dotsb x_n^{\alpha_{in}}x^{\beta_j} \notag \\
+ &\ m_i x_1^{\alpha_{i1}}p_{\alpha_{i2}, \sigma_3^{\alpha_{i3}}(\dotsb (\sigma_{{in}}^{\alpha_{in}}(b)))} x_3^{\alpha_{i3}}\dotsb x_n^{\alpha_{in}}x^{\beta_j} \notag \\
+ &\ m_i x_1^{\alpha_{i1}}x_2^{\alpha_{i2}}p_{\alpha_{i3}, \sigma_{i4}^{\alpha_{i4}} (\dotsb (\sigma_{in}^{\alpha_{in}}(b)))} x_4^{\alpha_{i4}}\dotsb x_n^{\alpha_{in}}x^{\beta_j}\notag  \\
+ &\ \dotsb + _i x_1^{\alpha_{i1}}x_2^{\alpha_{i2}} \dotsb x_{i(n-2)}^{\alpha_{i(n-2)}}p_{\alpha_{i(n-1)}, \sigma_{in}^{\alpha_{in}}(b)}x_n^{\alpha_{in}}x^{\beta_j} \notag \\
+ &\ m_i x_1^{\alpha_{i1}}\dotsb x_{i(n-1)}^{\alpha_{i(n-1)}}p_{\alpha_{in}, b}x^{\beta_j}.\label{perooo}
\end{align}
This guarantees that $M\langle X\rangle$ is really an $A$-module. In this way, when we compute every summand of $m_ix_1^{\alpha_{i1}}\dotsb x_n^{\alpha_{in}}b_jx_1^{\beta_{j1}}\dotsb x_n^{\beta_{jn}}$, we obtain products of the coefficient $m_i$ with several evaluations of $b_j$ in $\sigma$'s and $\delta$'s, depending of the coordinates of $\alpha_i$.\\

The purpose in the next two sections is to study the existing relations between an $R$-module $M_R$ and the polynomial module $M\langle X\rangle$ over the skew-PBW extension $A$ of $R$. Therefore, we extend the notions of skew-Armendariz modules and skew quasi-Armendariz modules introduced by Alhevaz and Moussavi \cite{AlhevazMoussavi2012} for the case of Ore extensions, and hence we generalize the concepts of $\alpha$-skew Armendariz modules \cite{ZhangChen2008} and $\alpha$-reduced modules \cite{LeeZhou2004}.
\subsection{Skew-Armendariz modules}\label{skewarmendarizmodules}
In this section we introduce the notion of skew-Armendariz module for $\sigma$-PBW extensions. As we said above, our treatment generalize \cite{LeeZhou2004}, \cite{ZhangChen2008}, and \cite{AlhevazMoussavi2012}. \\

Let us briefly recall some definitions about the notion of Armendariz for modules: (i) (\cite{ZhangChen2008}, Definition 2.1); let $B$ be a ring with an endomorphism $\alpha$ and $M_B$ an $B$-module. $M_B$ is called an $\alpha$-{\em skew Armendariz module}, if for polynomials $m(x)=m_0 + m_1x + \dotsb + m_kx^{x}\in M[X]$ and $f(x) = b_0 + b_1x + \dotsb + b_nx^{n}\in B[x;\alpha]$, $m(x)f(x) = 0$ implies $m_i\alpha^{i}(b_j) = 0$, for every $0\le i\le k$ and $0\le j\le n$. (ii) (\cite{AlhevazMoussavi2012}, Definition 2.2); let $B$ be a ring with an endomorphism $\alpha$ and $\alpha$-derivation $\delta$. Let $M_B$ be an $B$-module. $M_B$ is an $(\alpha, \delta)$-{\em skew Armendariz module}, if for polynomials $m(x) = m_0 + m_1x + \dotsb + m_kx^{k}\in M[X]$ and $f(x) = b_0 + b_1x+\dotsb + b_nx^{n} \in B[x;\alpha, \delta]$, $m(x)f(x) = 0$ implies $m_ix^{i}b_jx^{j}=0$, for every $0\le i\le k$ and $0\le j\le n$ (it is clear that if $\delta=0$, this definition coincides with the definition of $\alpha$-skew Armendariz module). (iii) (\cite{AlhevazMoussavi2012}, Definition 2.3); let $B$ be a ring with an endomorphism $\alpha$ and a $\alpha$-derivation $\delta$. A right $B$-module $M_B$ is called {\em skew-Armendariz module}, if for polynomials $m(x) = m_0 + m_1x+\dotsb + m_kx^{k}\in M[X]$ and $f(x)=b_0 + b_1x + \dotsb + b_nx^{n}\in B[x;\alpha, \delta]$, the equality $m(x)f(x) = 0$ implies $m_0b_j = 0$, for every $0\le j\le n$. All these definitions are generalized in the next definition for the context of $\sigma$-PBW extensions.
\begin{definition}\label{Generaldefinition2.3}
Let $A$ be a $\sigma$-PBW extension of a ring $R$, and let $M_R$ be a right $R$-module. $M_R$ is called a {\em skew-Armendariz module}, if for elements $m=m_0 + m_1X_1 + \dotsb + m_kX_k\in M\langle X\rangle$ and $f=b_0+b_1Y_1 + \dotsb + b_tY_t\in A$ with  $mf=0$, we have $m_0b_j=0$, for every $0\le j\le t$.
\end{definition}
Note that $B$ is skew-Armendariz if $B_B$ is a skew-Armendariz module.  Now, since the notion of skew-Armendariz module is a generalization of an $\alpha$-skew-Armendariz module (\cite{AlhevazMoussavi2012}, Theorem 2.4), both concepts in the context of Ore extensions, and our notion of skew-Armendariz module in Definition \ref{Generaldefinition2.3} is formulated for $\sigma$-PBW extensions, which are more general than Ore extensions (with $\alpha$ injective), then our skew-Armendariz module notion is more general than $\alpha$-skew-Armendariz module. Nevertheless, we can establish the following result without proof. Theorem \ref{generAlhevaz2012Theorem2.4} generalizes \cite{AlhevazMoussavi2012}, Theorem 2.4, Corollaries 2.5, 2.6 and 2.7.
\begin{theorem}\label{generAlhevaz2012Theorem2.4}
If $A$ is a quasi-commutative skew PBW extension of a ring $R$, and $M_R$ is a right module, then $M_R$ is $\Sigma$-skew Armendariz if and only if for every polynomials $m=m_0 + m_1X + \dotsb + m_kX_k\in M\langle X\rangle$, and $f=b_0 + b_1X_1+\dotsb + b_mX_m\in A$, the equality $mf=0$ implies $m_0b_j=0$, for every $0\le j\le m$.
\end{theorem}
The next definition introduce a more general class of modules than those established in Definition \ref{Generaldefinition2.3}.
\begin{definition}\label{Generaldefinition2.4}
Let $A$ be a $\sigma$-PBW extension of a ring $R$, and let $M_R$ be a right $R$-module. $M_R$ is called a {\em linearly skew-Armendariz module}, if for linear polynomials $m=m_0 + m_1x + \dotsb + m_nx_n\in M\langle X\rangle$ and $g(x)=b_0 + b_1x + \dotsb + b_nx_n\in A$ with $mg=0$, we have $m_0b_j=0$, for every $0\le j\le n$.
\end{definition}
In \cite{AnninPhD2002}, Annin introduce the notion of compatibility for modules in the following way: given a module $M_B$, an endomorphism $\alpha:B\to B$ and an $\alpha$-derivation $\delta:R\to R$, $M_B$ is $\alpha$-compatible if for each $m\in M$ and $r\in B$, we have $mr=0\Leftrightarrow m\alpha(r)=0$. Moreover, $M_R$ is $\delta$-compatible if for each $m\in M$ and $r\in B$, we have $mr=0 \Rightarrow m\delta(r)=0$. If $M_B$ is both $\alpha$-compatible and $\delta$-compatible, $M_B$ is called $(\alpha, \delta)$-compatible. In \cite{ReyesSuarezUMA2017}, Definition 3.2, the author defined the notion of compatibility for skew PBW extensions in the following way (this definition extends \cite{HashemiMoussavi2005}): consider a ring $R$ with a family of endomorphisms  $\Sigma$ and a family of $\Sigma$-derivations $\Delta$ (Proposition \ref{sigmadefinition}). (i)  $R$ is said to be $\Sigma$-{\em compatible}, if for each $a,b\in R$, $a\sigma^{\alpha}(b)=0$ if and only if $ab=0$, for every $\alpha\in \mathbb{N}^{n}$; (ii) $R$ is said to be $\Delta$-{\em compatible}, if for each $a,b \in R$, $ab=0$ implies $a\delta^{\beta}(b)=0$, for every $\beta \in \mathbb{N}^{n}$; (iii)
if $R$ is both $\Sigma$-compatible and $\Delta$-compatible, $R$ is called $(\Sigma, \Delta)$-{\em compatible}. As it was established in \cite{ReyesSuarezUMA2017}, Proposition 3.3, the importance of $(\Sigma, \Delta)$-compatible rings is that they are more general than $\Sigma$-rigid rings defined and characterized by the author in \cite{Reyes2015} in terms of the properties of being Baer, p.p., p.q., and p.q.-Baer ($\Sigma$-rigid rings are a generalization of $\alpha$-rigid rings defined by Krempa in \cite{Krempa1996} and studied by Hong et. al., \cite{HongKimKwak2000}). Next, we extend this definition of compatibility for the context of modules over skew PBW extensions.
\begin{definition}
If $A$ is a $\sigma$-PBW extension of a ring $R$, and $M_R$ is a right $R$-module, then $M_R$ is called $\Sigma$-{\em compatible}, if for every $m\in M$ and $r\in R$, $mr=0\Leftrightarrow m\sigma^{\alpha}(r)=0$, for any $\alpha \in \mathbb{N}^{n}$. $M_R$ is called $\Delta$-{\em compatible}, if for every $m\in M$ and $r\in R$, $mr=0 \Rightarrow m\delta^{\beta}(r)=0$, for any $\beta\in \mathbb{N}^{n}$. If $M_R$ is both $\Sigma$-compatible and $\Delta$-compatible, then $M_R$ is called $(\Sigma, \Delta)$-{\em compatible}.
\end{definition}
From \cite{AnninPhD2002}, Lemma 2.16, we know that in the case of Ore extensions, a module $M_B$ is $(\alpha, \delta)$-compatible if and only if the polynomial extension $M\langle X\rangle_B$ is $(\alpha, \delta)$-compatible. This assertion os generalized in the next proposition.
\begin{proposition}
If $A$ is a $\sigma$-PBW extension of $R$, and $M_R$ is a right $R$-module, then $M_R$ is $(\Sigma, \Delta)$-compatible if and only if $M\langle X\rangle_R$ is $(\Sigma, \Delta)$-compatible.
\begin{proof}
The proofs follow from the definitions.
\end{proof}
\end{proposition}
The following proposition is a direct consequence from \cite{AlhevazMoussavi2012}, Lemma 2.14.
\begin{proposition}\label{generAlhevazMoussavi2012Lemma2.14}
If $M_R$ is a right $R$-module, then the following conditions are equivalent:
\begin{enumerate}
\item [\rm (i)] $M_R$ is reduced and $(\Sigma, \Delta)$-compatible.
\item [\rm (ii)] for any $m\in M$ and $r\in R$, the following conditions hold:
\begin{enumerate}
\item [\rm (a)] $mr=0$ implies $mRr=0$;
\item [\rm (b)] $mr=0$ implies $m\delta^{\beta}(r)=0$, for any $\beta\in \mathbb{N}^{n}$;
\item [\rm (c)] $mr=0$ if and only if $m\sigma^{\theta}(r)=0$, for any $\theta\in \mathbb{N}^{n}$;
\item [\rm (d)] $mr^{2} = 0$ implies $mr=0$.
\end{enumerate}
\end{enumerate}
\end{proposition}
The following proposition generalizes \cite{AlhevazMoussavi2012}, Lemma 2.15.
\begin{proposition}\label{generAlehevaz2012Lemma2.15}
If $M_R$ is an $(\Sigma, \Delta)$-compatible module, and $m\in M,\ a, b\in R$, then we have the following assertions:
\begin{enumerate}
\item [\rm (i)] if $ma=0$, then $m\sigma^{\theta}(a) = 0 = m\delta^{\theta}(a)$, for any element $\theta\in \mathbb{N}^{n}$;
\item [\rm (ii)] if $mab=0$, then $m\sigma_i(\delta^{\theta}(a))\delta_i(b)=m\sigma^{\beta}(\delta_i(a))\delta^{\theta}(b)$, and so, $ma\delta^{\theta}(b) = 0 = m\delta^{\theta}(a)b$, for any elements $\beta, \theta\in \mathbb{N}^{n}$, and $i=1,\dotsc, n$;
\item [\rm (iii)] ${\rm ann}_R(\{ma\}) = {\rm ann}_R(m\sigma_i(a)) = {\rm ann}_R(\{m\delta_i(a)\})$, for every $i=1,\dotsc, n$.
\end{enumerate}
\end{proposition}
Proposition \ref{generAlhevaz2012Lemma2.16} generalizes  \cite{AlhevazMoussavi2012}, Lemma 2.16, from Ore extensions to $\sigma$-PBW extensions.
\begin{proposition}
\label{generAlhevaz2012Lemma2.16}
If $A$ is a $\sigma$-PBW extension of a ring $R$, $M_R$ is a $(\Sigma, \Delta)$-compatible right $R$-module,  $m=m_0 + m_1X_1 + \dotsb + m_kX_k$ is an element of $M\langle X\rangle$, and $a\in B$, then $mr=0$ if and only if $m_ir=0$, for every $0\le i\le k$.
\begin{proof}
Suppose that $m_ir=0$, for every $0\le i\le k$. Since
\begin{align}
mr = &\ (m_0 + m_1X_1 + \dotsb + m_kX_k)r\notag\\ =&\ m_0r +
m_1X_1r + \dotsb + m_kX_kr\notag\\ =&\ m_0r +
m_1(\sigma^{\alpha_1}(r)X_1 + p_{\alpha_1, r}) + \dotsb +
m_k(\sigma^{\alpha_k}(r)X_k + p_{\alpha_k, r})\notag\\ =&\ m_0r +
m_1\sigma^{\alpha_1}(r)X_1 + m_1p_{\alpha_1,r} + \dotsb +
m_k\sigma^{\alpha_k}(r)X_k + m_kp_{\alpha_k,r},\label{qwerm}
\end{align}
where $\alpha_i={\rm exp}(X_i)$, $p_{\alpha_i,r} = 0$, or, ${\rm deg}(p_{\alpha_i,r}) < |\alpha_i|$, if $p_{\alpha_i,r}\neq 0$, for every $i$, and using the equality $m_ir=0$ with the expression (\ref{perooo}) and the $(\Sigma, \Delta)$-compatibility of $M_R$, we conclude that $mr=0$.

Now, suppose that $mr=0$. From expression (\ref{qwerm}) we can see that ${\rm lc}(mr)=m_k\sigma^{\alpha_k}(r)$, so by the $\Sigma$-compatibility of $M_R$, we obtain $m_kr=0$. Hence, expresion (\ref{perooo}) and $(\Sigma,\Delta)$-compatibility of $M_R$ imply that $p_{\alpha_k,r}=0$, so $mr$ reduces to
\[
mr = m_0r + m_1\sigma^{\alpha_1}(r)X_1 + m_1p_{\alpha_1,r} + \dotsb + m_{k-1}\sigma^{\alpha_{k-1}}(r)X_{k-1} + m_{k-1}p_{\alpha_{k-1},r}.
\]
Again, since ${\rm lc}(mr) = m_{k-1}\sigma^{\alpha_{k-1}}(r)=0$, from $\Sigma$-compatibility of $M_R$ we can assert that $m_{k-1}r=0$. In this way, expresion (\ref{perooo}) and $(\Sigma,\Delta)$-compatibility of $M_R$ imply that $p_{\alpha_{k-1},r}=0$, so $mr$ takes the form
\[
mr = m_0r + m_1\sigma^{\alpha_1}(r)X_1 + m_1p_{\alpha_1,r} + \dotsb + m_{k-2}\sigma^{\alpha_{k-2}}(r)X_{k-2} + m_{k-2}p_{\alpha_{k-2},r}.
\]
Continuing in this way we can show  that $m_kr = m_{k-1}r = m_{k-2}r =  \dotsb = m_1r = m_0r$, which concludes the proof.
\end{proof}
\end{proposition}
The next proposition generalizes \cite{AlhevazMoussavi2012}, Proposition 2.17 and Corollary 2.18.
\begin{proposition}\label{generAlhevaz2012Proposition2.17}
A module $M_R$ is $\Sigma$-reduced if and only if the extension $M\langle X\rangle_R$ is an $\Sigma$-reduced module.
\end{proposition}
Theorem \ref{generAlhevaz2012Theorem2.19} generalizes  \cite{AlhevazMoussavi2012}, Theorem 2.19.
\begin{theorem}\label{generAlhevaz2012Theorem2.19}
If $M_R$ is an $(\Sigma, \Delta)$-compatible and reduced module, then $M_R$ is skew-Armendariz.
\begin{proof}
Consider the elements $m=m_0 + m_1X_1 + \dotsb + m_kX_k\in M\langle X\rangle,\ f=b_0 + b_1Y_1 + \dotsb + a_tY_t\in A$, with $mf=0$. We have $mf= (m_0+m_1X_1+\dotsb + m_kX_k)(b_0+b_1Y_1+\dotsb + b_tY_t)=\sum_{l=0}^{k+t} \biggl(\sum_{i+j=l} m_iX_ib_jY_j\biggr)$. Note that ${\rm lc}(mf)= m_k\sigma^{\alpha_k}(b_t)c_{\alpha_k, \beta_t}=0$. Since $A$ is bijective, $m_k\sigma^{\alpha_k}(b_t)=0$, and by the $\Sigma$-compatibility of $M_R$, $m_kb_t=0$. The idea is to prove that $m_pb_q=0$ for $p+q\ge 0$.  We proceed  by induction. Suppose that $m_pb_q=0$ for $p+q=k+t, k+t-1, k+t-2, \dotsc, l+1$ for some $l>0$. By Proposition \ref{generAlhevazMoussavi2012Lemma2.14} and expression (\ref{perooo}), we obtain $m_pX_pb_qY_q=0$, for these values of $p+q$. In this way, we only consider the sum of the products $m_uX_ub_vY_v$, where $u+v=l, l-1,l-2,\dotsc, 0$. Fix $u$ and $v$. Consider the sum of all terms of $mf$  having exponent $\alpha_u+\beta_v$. From expression (\ref{perooo}), Proposition \ref{generAlhevazMoussavi2012Lemma2.14}, and the assumption $mf=0$, we know that the sum of all coefficients of all these terms  can be written as
\begin{equation}\label{Federer}
m_u\sigma^{\alpha_u}(b_v)c_{\alpha_u, \beta_v} + \sum_{\alpha_{u'} + \beta_{v'} = \alpha_u + \beta_v} m_{u'}\sigma^{\alpha_{u'}} ({\rm \sigma's\ and\ \delta's\ evaluated\ in}\ b_{v'})c_{\alpha_{u'}, \beta_{v'}} = 0.
\end{equation}
By assumption, we know that $m_pb_q=0$, for $p+q=k+t, k+t-1, \dotsc, l+1$.  So, Proposition \ref{generAlhevazMoussavi2012Lemma2.14} guarantees that the product
\[m_p({\rm \sigma's\ and\ \delta's\ evaluated\ in}\ b_{q})\ \ \ \ \ \ \ ({\rm any\  order\ of}\ \sigma's\ {\rm and}\ \delta's)
\]
is equal to zero. Then  $[({\rm \sigma's\ and\ \delta's\ evaluated\ in}\ b_{q})a_p]^2=0$, and hence we obtain the equality $({\rm \sigma's\ and\ \delta's\ evaluated\ in}\ b_{q})m_p=0$ ($M_R$ is reduced). In this way, multiplying (\ref{Federer}) by $m_l$, and using the fact that the elements $c_{i,j}$ in Definition \ref{gpbwextension} (iv) are in the center of $R$,
\begin{equation}\label{doooooo}
m_u\sigma^{\alpha_u}(b_v)a_kc_{\alpha_u, \beta_v} + \sum_{\alpha_{u'} + \beta_{v'} = \alpha_u + \beta_v} m_{u'}\sigma^{\alpha_{u'}} ({\rm \sigma's\ and\ \delta's\ evaluated\ in}\ b_{v'})m_lc_{\alpha_{u'}, \beta_{v'}} = 0,
\end{equation}
whence, $m_u\sigma^{\alpha_u}(b_0)a_l=0$. Since $u+v=l$ and $v=0$, then $u=l$, so $m_l\sigma^{\alpha_l}(b_0)m_l=0$, i.e., $[m_l\sigma^{\alpha_l}(b_0)]^{2}=0$, from which $m_l\sigma^{\alpha_l}(b_0)=0$ and $m_lb_0=0$, by Proposition \ref{generAlhevazMoussavi2012Lemma2.14}.  Therefore, we now have to study the expression (\ref{Federer}) for $0\le u \le l-1$ and $u+v=l$. If we multiply (\ref{doooooo}) by $m_{l-1}$, we obtain
{\small{\[
m_u\sigma^{\alpha_u}(b_v)m_{l-1}c_{\alpha_u, \beta_v} + \sum_{\alpha_{u'} + \beta_{v'} = \alpha_u + \beta_v} m_{u'}\sigma^{\alpha_{u'}} ({\rm \sigma's\ and\ \delta's\ evaluated\ in}\ b_{v'})m_{k-1}c_{\alpha_{u'}, \beta_{v'}} = 0.
\]}}
Using a similar reasoning as above, we can see that $m_u\sigma^{\alpha_u}(b_1)m_{l-1}c_{\alpha_u, \beta_1}=0$. Since $A$ is bijective, $m_u\sigma^{\alpha_u}(b_1)m_{l-1}=0$, and using the fact $u=l-1$, we have $[m_{l-1}\sigma^{\alpha_{l-1}}(b_1)]=0$, which imply $m_{l-1}\sigma^{\alpha_{l-1}}(b_1)=0$, that is, $m_{l-1}b_1=0$. Continuing in this way we prove that $m_ib_j=0$, for $i+j=l$. Therefore $a_ib_j=0$, for $0 \le i\le m$ and $0\le j\le t$, which concludes the proof.
\end{proof}
\end{theorem}
\begin{proposition}[\cite{AlhevazMoussavi2012}, Proposition 2.22]
Suppose that $M$ is a flat right $B$-module. Then for every exact sequence $0\to K\to F\to M\to 0$, where $F$ is $B$-free, we have $FI\cap K = KI$ for each left ideal $I$ of $B$. In particular, we have $Fa\cap K = Ka$ for each element $a$ of $B$.
\end{proposition}
Proposition \ref{generAlhevaz2012Proposition2.23} generalizes
\cite{AlhevazMoussavi2012}, Proposition 2.23.
\begin{proposition}\label{generAlhevaz2012Proposition2.23}
If $A$ is a $\sigma$-PBW extension of a ring $R$, then $B$ is a skew-Armendariz ring if and only if every flat $B$-module $M$ is skew-Armendariz.
\begin{proof}
Consider $M$ a flat $R$-module, and an exact sequence $0\to K\to F\to M\to 0$, where $F$ is free over $R$. If $b$ is an element of $F$, then $\overline{b} = b + K$ is an element of $M$. Let $f = \overline{b}_0 + \overline{b_1}X_1 + \dotsb + \overline{b_k}X_k\in M\langle X\rangle$ and $g = a_0 + a_1Y_1 + \dotsb + a_tY_t\in A$ with $fg=0$. Let us prove that $\overline{b}_0a_j = 0$, for $0\le j\le m$. From the expression for the product $fg$ given by $fg = (\overline{b}_0 + \overline{b_1}X_1 + \dotsb + \overline{b_t}X_t) (a_0 + a_1X_1 + \dotsb + a_mX_m) = \sum_{l=0}^{k+t} \biggl(\sum_{i+j=l} \overline{b}_iX_ia_jY_j\biggr) $, and using the relations established in (\ref{perooo}), we can find explicitly the coefficients for every term of $fg$ (considering a total order as in Definition \ref{definitioncoefficients} for the products of the elements $X_i,\ Y_j$, for instance, $X_k\succ \dotsb \succ X_1$ and $Y_t \succ \dotsb \succ Y_1$). Now, by assumption, $M$ is a flat $R$-module, so there exists an $R$-module homomorphism $\beta: F\to K$ which fixes the coefficients of every term of the product $fg$. Consider the elements $w_i:=\beta(b_i) - b_i,\ i=1,\dotsc, k$ in $F$. Then the element $h = w_0 + w_1X_1 + \dotsb + w_kX_k$ is an element of $F[X]$ which satisfies that $hg=0$. Note that $F$ is skew-Armendariz by Proposition, because $R$ is skew-Armendariz and $F$ is a free $R$-module. Hence, $w_0a_j = 0$, for $j=1,\dotsc, k$, and so $b_0a_j\in K$, for every $j$, that is, $\overline{b}_0a_j=0$ in $M$, which shows that $M$ is skew-Armendariz.
\end{proof}
\end{proposition}
For the next theorem, consider the set ${\rm Ann}_B(2^{M_B}) := \{{\rm ann}_B(U)\mid U\subseteq M_B\}$, where $M_B$ is an $B$-module. Theorem \ref{generAlhevazTheorem 2.24} generalizes \cite{AlhevazMoussavi2012}, Theorem 2.24, from Ore extensions to skew PBW extensions.
\begin{theorem}\label{generAlhevazTheorem 2.24}
If $A$ is a $\sigma$-PBW extension of a ring $R$, and $M_R$ is a $(\Sigma, \Delta)$-compatible right $R$-module, then the following conditions are equivalent:
\begin{enumerate}
\item [\rm (1)] $M_R$ is a skew-Armendariz module;
\item [\rm (2)] The map $\psi:{\rm Ann}_R(2^{M_R})\to {\rm Ann}_S(2^{M\langle X\rangle_A})$, defined by $C\to CA$, for all $C\in {\rm Ann}_R(2^{M_R})$, is bijective.
\end{enumerate}
\begin{proof}
$(1)\Rightarrow (2)$ First of all, let us see that ${\rm ann}_R(U)A = {\rm ann}_A(U)$, for every $U\subseteq M_R$. If $f\in {\rm ann}_R(U)A$, then $f$ is expressed as $f=r_0 + r_1Y_1 + \dotsb + r_kY_p$, where $r_i\in {\rm ann}_R(U)$ and $Y_i\in {\rm Mon}(A)$, for every $i$. If $g = m_0 + m_1X_1 + \dotsb + m_kX_k\in U$, then $gf = 0$ because $gr_i = 0$, for $1\le i\le p$, so $gf=0$, i.e., $f\in {\rm ann}_A(U)$. Now, if $h = a_0 + a_1Y_1 + \dotsb + a_tY_t\in {\rm ann}_A(U)$, then $gh=0$, for every $g = m_0 + m_1X_1 + \dotsb + m_kX_k\in U$. Since $M_R$ is a skew-Armendariz module, we have $m_0a_j=0$, for $j=0,\dotsc, t$, and by the $(\Sigma, \Delta)$-compatibility of $M_R$, we can see that $ga_j = 0$, for every $j$, so $h\in {\rm ann}_R(U)A$.

The application $\psi:\{{\rm ann}_R(U)\mid U\subseteq M_R\}\to \{{\rm ann}_A(U)\mid U\subseteq M\langle X\rangle_A\}$, defined by $C\to CA$, for every $C\in \{{\rm ann}_R(U)\mid U\subseteq M_R\}$, is well defined,  since ${\rm ann}_R(U)A = {\rm ann}_A(U)$, for every $U\subseteq M_R$. From the $(\Sigma, \Delta)$-compatibility of $M_R$, we can deduce that ${\rm ann}_A(V)\cap R = {\rm ann}_R(V_0)$, for every $V\subseteq M\langle X\rangle_A$, where $V_0\subseteq M$ is the set of coefficients of $V$. This fact guarantees that the map $\psi':\{{\rm ann}_A(U)\mid U\subseteq M\langle X\rangle_A\}\to \{{\rm ann}_R(U)\mid U\subseteq M_R\}$, defined by $D\to D\cap R$, for every $D\in \{{\rm ann}_A(U)\mid U\subseteq M\langle X\rangle_A\}$. Note that $\psi'\psi={\rm id}$, so $\psi$ is an injective map. Consider the set $B\in \{{\rm ann}_A(U)\mid U\subseteq M\langle X\rangle_A\}$. Then $B={\rm ann}_A(J)$, for some $J\subseteq M\langle X\rangle_A$. Let $B_1$ and $J_1$ be the set of coefficients of elements of $B$ and $J$, respectively. The aim is to prove that ${\rm ann}_R(J_1)=B_1R$. With this in mind, let $m = m_0 + m_1X_1 + \dotsb + m_kX_k\in J$ and $f=a_0 + a_1Y_1 + \dotsb + a_tY_t\in B$. Then $mf=0$, and using both  assumptions on $M_R$, it follows that $m_ib_j=0$, for every $m_i$ and $b_j$, whence $J_1B_1 = 0$, that is, $B_1R\subseteq {\rm ann}_R(J_1)$. Note that due to the $(\Sigma, \Delta)$-compatibility of $M_R$, we can assert that ${\rm ann}_R(J_1)\subseteq B_1R$, which shows that ${\rm ann}_R(J_1)=B_1R$, and so ${\rm ann}_A(J) = B_1RA$. This proves that $\psi$ is a surjective function. 

$(2) \Rightarrow (1)$ Consider the elements $m = m_0 + m_1X_1 + \dotsb + m_kX_k\in M\langle X\rangle_A$ and $f\in a_0 + a_1Y_1 + \dotsb + a_tY_t\in A$ with $mf=0$. It is clear that $f\in {\rm ann}_A(\{m\})={\rm ann}_R(U)A$, where $U\subseteq M_R$. In this way, the elements $b_0,\dotsc, b_n$ belong to ${\rm ann}_R(U)$, whence $mb_j=0$, for every $j$. Therefore, $m_0b_j=0$, for $0\le j\le t$, which concludes the proof.
\end{proof}
\end{theorem}
The following theorem generalizes \cite{AlhevazMoussavi2012}, Theorem 2.25.
\begin{theorem}\label{generAlhevaz2012Theorem2.25}
If $A$ is a $\sigma$-PBW extension of a ring  $R$, and $M_R$ is a linearly skew-Armendariz module with $R\subseteq M$, then for every idempotent elements $e$ of $R$, we have $\sigma_i(e)=e$ and $\delta_i(e)=0$, for every $i=1,\dotsc, n$.
\begin{proof}
Note that $M_R$ is a linearly skew-Armendariz module with $R\subseteq M_R$, so $R_R$ is also linearly skew-Armendariz. Let us prove the assertion for $R_R$.

Consider an idempotent element $e$ of $R$. Then $\delta_i(e)=\sigma_i(e)\delta_i(e)+\delta_i(e)e$. Let $f, g\in A$ given by $f = \delta_i(e) + 0x_1 + \dotsb + 0x_{i-1} + \sigma_i(e)x_i + 0x_{i+1} + \dotsb + 0x_n$, and $g = e-1 + (e-1)x_1 + \dotsb + (e-1)x_n$, respectively.
Recall that $\delta_i(1)=0$, for every $i$. Let us show that $fg=0$:
\begin{align*}
fg = &\ \delta_i(e)(e-1) + \biggl(\sum_{j=1}^{n} \delta_i(e)(e-1)x_j\biggr) + \sigma_i(e)x_i(e-1) + \sum_{j=1}^{n} \sigma_i(e)x_i(e-1)x_j\\
= &\ \delta_i(e)(e-1) + \biggl(\sum_{j=1}^{n} \delta_i(e)(e-1)x_j\biggr) + \sigma_i(e)[\sigma_i(e-1)x_i + \delta_i(e-1)] \\
+ &\ \sum_{j=1}^{n} \sigma_i(e)[\sigma_i(e-1)x_i + \delta_i(e-1)]x_j.
\end{align*}
Equivalently,

\begin{align*}
fg = &\ \delta_i(e)(e-1) + \biggl(\sum_{j=1}^{n} \delta_i(e)(e-1)x_j\biggr) + \sigma_i(e) [(\sigma_i(e) - \sigma_i(1))x_i + \delta_i(e)]\\
+ &\ \sum_{j=1}^{n} \sigma_i(e) [(\sigma_i(e) - \sigma_i(1))x_i + \delta_i(e)]x_j\\
= &\ \delta_i(e)(e-1) + \biggl(\sum_{j=1}^{n} \delta_i(e)(e-1)x_j\biggr) + \sigma_i(e)[\sigma_i(e)x_i - x_i +\delta_i(e)]\\
+ &\ \sum_{j=1}^{n} \sigma_i(e)[\sigma_i(e)x_i - x_i +\delta_i(e)]x_j\\
= &\ \delta_i(e)e - \delta_i(e) + \biggl(\sum_{j=1}^{n}(\delta_i(e)e - \delta_i(e))x_j\biggr) + \sigma_i(e)x_i - \sigma_i(e)x_i + \sigma_i(e)\delta_i(e)\\
+ &\ \sum_{j=1}^{n} (\sigma_i(e)x_i - \sigma_i(e)x_i + \sigma_i(e)\delta_i(e))x_j\\
= &\ \delta_i(e)e - \delta_i(e) + \sum_{j=1}^{n} \delta_i(e)ex_j - \sum_{j=1}^{n} \delta_i(e)x_j + \sigma_i(e)\delta_i(e) + \sum_{j=1}^{n} \sigma_i(e)\delta_i(e)x_j = 0.
\end{align*}
From Definition \ref{Generaldefinition2.4} we obtain $\delta_i(e)(e-1) = 0$, i.e., $\delta_i(e)e = \delta_i(e)$, and hence $\sigma_i(e)\delta_i(e) = 0$.

Now, consider the elements $s$ and $t$ of $A$ given by $s=  \delta_i(e) - (1-\sigma_i(e))x_i$ and $t=e + \sum_{j=1}^{n} ex_j$, respectively. Let us show that $st=0$:
\begin{align*}
st = &\ \delta_i(e)e + \biggl(\delta_i(e)e\sum_{j=1}^{n} x_j\biggr) - (1-\sigma_i(e))x_ie - \biggl((1-\sigma_i(e))x_ie\sum_{j=1}^{n} x_j\biggr)\\
= &\ \delta_i(e)e + \biggl(\delta_i(e)e\sum_{j=1}^{n}x_j\biggr) - x_ie + \sigma_i(e)x_ie - x_ie\sum_{j=1}^{n} x_j + \sigma_i(e)x_ie\sum_{j=1}^{n} x_j\\
= &\ \delta_i(e)e + \biggl(\delta_i(e)e\sum_{j=1}x_j\biggr) - (\sigma_i(e)x_i + \delta_i(e)) + \sigma_i(e)(\sigma_i(e)x_i + \delta_i(e))\\
- &\ \biggl((\sigma_i(e)x_i + \delta_i(e))\sum_{j=1}^{n} x_j\biggr) + \sigma_i(e)(\sigma_i(e)x_i + \delta_i(e))\sum_{j=1}^{n} x_j\\
= &\ \delta_i(e)e + \biggl(\delta_i(e)e\sum_{j=1}^{n}x_j\biggr) - \sigma_i(e)x_i - \delta_i(e) + \sigma_i(e)x_i + \sigma_i(e)\delta_i(e) - \sigma_i(e)x_i\sum_{j=1}^{n}x_j \\
- &\ \delta_i(e)\sum_{j=1}^{n}x_j + \sigma_i(e)x_i\sum_{j=1}^{n} x_j + \sigma_i(e)\delta_i(e)\sum_{j=1}^{n}x_j.
\end{align*}
Since $\delta_i(e)=\delta_i(e)e$ and $\sigma_i(e)\delta_i(e)=0$, then $st=0$. By Armendariz condition we know that $\delta_i(e)e=0$, which shows that $\delta_i(e)=0$.

Consider the elements $u, v \in A$ given by $u=1-e + (1-e)\sigma_i(e)x_i$ and $v=e+(e-1)\sigma_i(e)x_i$. We have the equalities
\begin{align*}
uv = &\ e+ (e-1)\sigma_i(e)x_i - e^2 - e(e-1)\sigma_i(e)x_i + (1-e)\sigma_i(e)x_ie + (1-e)\sigma_i(e)x_i(e-1)\sigma_i(e)x_i\\
= &\ e\sigma_i(e)x_i - \sigma_i(e)x_i - e\sigma_i(e)x_i + e\sigma_i(e)x_i + (1-e)\sigma_i(e)(\sigma_i(e)x_i + \delta_i(e))\\
+ &\ (1-e)\sigma_i(e)(\sigma_i(e)x_i - x_i + \delta_i(e))\sigma_i(e)x_i\\
= &\ -\sigma_i(e)x_i + e\sigma_i(e)x_i + \sigma_i(e)x_i + \sigma_i(e)\delta_i(e) - e\sigma_i(e)x_i - e\sigma_i(e)\delta_i(e) \\
+ &\ [\sigma_i(e)x_i - \sigma_i(e)x_i + \sigma_i(e)\delta_i(e) - e\sigma_i(e)x_i + e\sigma_i(e)x_i - e\sigma_i(e)\delta_i(e)]\sigma_i(e)x_i = 0.
\end{align*}
Using that $\delta_i(e)=0$, we obtain $(1-e)(e-1)\sigma_i(e)=0$, i.e., $e\sigma_i(e) = \sigma_i(e)$.

Finally, let $w=e+e(1-\sigma_i(e))x_i,\ z=1-e - e(1-\sigma_i(e))x_i$ be elements of $A$. Then
\begin{align*}
wz = &\ e-e^2-e^2(1-\sigma_i(e))x_i + e(1-\sigma_i(e))x_i - e(1-\sigma_i(e))x_ie - e(1-\sigma_i(e))x_ie(1-\sigma_i(e))x_i\\
= &\ -e(1-\sigma_i(e))(\sigma_i(e)x_i + \delta_i(e)) - e(1-\sigma_i(e))[\sigma_i(e(1-\sigma_i(e)))x_i + \delta_i(e(1-\sigma_i(e)))]x_i.
\end{align*}
Using that $\delta_i(e)=0$ and $e\sigma_i(e)=\sigma_i(e)$, we can see that $wz=0$. Hence, $e(-e(1-\sigma_i(e))) = 0$, which shows that $e\sigma_i(e)=e$, and so $\sigma_i(e)=e$.
\end{proof}
\end{theorem}
In  \cite{AgayevGungorogluHarmanciHalicioglu2009}, Agayev et. al., introduced the notion of abelian module: a right  $B$-module $M_B$ is called {\em abelian}, if for any elements $m\in M$, $r\in B$, and every idempotent $e\in B$, we have $mae=mea$. They proved that every Armendariz module, and hence every reduced module is abelian.
The next theorem generalizes  \cite{AlhevazMoussavi2012}, Theorem 2.26, from Ore Extensions to skew PBW extensions.
\begin{theorem}\label{generAlhevaz2012Theorem2.26}
If $A$ is a $\sigma$-PBW extension of a ring $R$, and $M_R$ is a linearly skew-Armendariz module with $R\subseteq M$, then $M_R$ is an abelian module.
\begin{proof}
Consider $M_R$ a linearly skew-Armendariz module, the elements defined by $m_1 = me - \sum_{i=1}^{n} mer(1-e)x_i,\ m_2=m(1-e) - \sum_{i=1}^{n}m(1-e)rex_i\in M\langle X\rangle_A$, and $f_1 = (1-e) + \sum_{i=1}^{n}er(1-e)x_i,\ f_2 = e + \sum_{i=1}^{n}(1-e)rex_i$, where $e\in R$ is an idempotent element, and $r\in R, m\in M$. Let us show that $m_1f_1=m_2f_2=0$. Recall that $\sigma_i(e)=e$ and $\delta_i(e)=0$, for $i=1,\dotsc, n$ (Theorem \ref{generAlhevaz2012Theorem2.25}), so $x_ie = ex_i$, for every $i$. We have the following equalities:
\begin{align*}
m_1f_1 = &\ \biggl(me - \sum_{i=1}^{n} mer(1-e)x_i \biggr) \biggl((1-e) + \sum_{i=1}^{n}er(1-e)x_i \biggr)\\
= &\ \biggl(me - \sum_{i=1}^{n} merx_i +\sum_{i=1}^{n} merex_i \biggr) \biggl(1-e + \sum_{i=1}^{n} erx_i - \sum_{i=1}^{n} erex_i \biggr)\\
= &\ me - me^{2} + \sum_{i=1}^{n}me^{2}rx_i - \sum_{i=1}^{n} me^{2}rex_i - \sum_{i=1}^{n} merx_i + \sum_{i=1}^{n} merx_ie \\
- &\ \biggl(\sum_{i=1}^{n} merx_i\biggr) \biggl(\sum_{i=1}^{n} erx_i\biggr) + \biggl(\sum_{i=1}^{n} merx_i\biggr)\biggl(\sum_{i=1}^{n} erex_i\biggr) - \sum_{i=1}^{n} merex_i + \sum_{i=1}^{n} merex_ie\\
+ &\ \biggl(\sum_{i=1}^{n}merex_i\biggr) \biggl(\sum_{i=1}^{n} erx_i\biggr) -  \biggl(\sum_{i=1}^{n} merex_i\biggr)\bigg(\sum_{i=1}^{n} erex_i\biggr),
\end{align*}
or equivalently,
\begin{align*}
m_1f_1= &\ - \biggl(\sum_{i=1}^{n}merex_ir\biggr) \biggl(\sum_{i=1}^{n}x_i\biggr) + \biggl(\sum_{i=1}^{n}merex_ir\biggr) \biggl(\sum_{i=1}^{n} ex_i\biggr)\\
+ &\ \biggl(\sum_{i=1}^{n} mere^{2}x_ir\biggr)\biggl(\sum_{i=1}^{n} x_i\biggr) - \biggl(\sum_{i=1}^{n} mere^{2}x_ir\biggr) \biggl(\sum_{i=1}^{n} ex_i\biggr) = 0,
\end{align*}
and,
\begin{align*}
m_2f_2 = &\ \biggl(m(1-e) - \sum_{i=1}^{n}m(1-e)rex_i\biggr) \biggl(e + \sum_{i=1}^{n}(1-e)rex_i \biggr)\\
= &\ \biggl(m-me - \sum_{i=1}^{n} mrex_i + \sum_{i=1}^{n}merex_i\biggr)\biggl(e + \sum_{i=1}^{n}rex_i - \sum_{i=1}^{n} erex_i\biggr)\\
= &\ me + \sum_{i=1}^{n}mrex_i - \sum_{i=1}^{n}merex_i - me^{2} - \sum_{i=1}^{n}merex_i + \sum_{i=1}^{n}me^{2}rex_i\\
- &\ \sum_{i=1}^{n}mrex_ie - \biggl(\sum_{i=1}^{n}mrex_i\biggr) \biggl(\sum_{i=1}^{n} rex_i \biggr) + \biggl(\sum_{i=1}^{n}mrex_i \biggr)\biggl(\sum_{i=1}^{n} erex_i\biggr)\\
+ &\ \sum_{i=1}^{n} merex_ie + \biggl(\sum_{i=1}^{n}merex_i\biggr)\biggl(\sum_{i=1}^{n}rex_i\biggr) - \biggl(\sum_{i=1}^{n}merex_i\biggr)\biggl(\sum_{i=1}^{n}erex_i\biggr)\\
= &\ - \biggl(\sum_{i=1}^{n}mrex_i\biggr)\biggl(\sum_{i=1}^{n} rex_i \biggr) + \biggl(\sum_{i=1}^{n}mrex_ie\biggr)\biggl(\sum_{i=1}^{n} rex_i \biggr)\\
+ &\ \biggl(\sum_{i=1}^{n}merex_i \biggr)\biggl(\sum_{i=1}^{n}
rex_i \biggr) - \biggl(\sum_{i=1}^{n}
merex_ie\biggr)\biggl(\sum_{i=1}^{n} rex_i\biggr) = 0.
\end{align*}
By assumption, $M_R$ is linearly skew-Armendariz, so $meer(1-e)=0$ and $m(1-e)(1-e)re$, or what is the same, $mer=mere$ and $mre=mere$, respectively, whence $mer=mre$, i.e., $M_R$ is an abelian module.
\end{proof}
\end{theorem}
\begin{corollary}\label{generAlhevaz2012Corollary2.27}
If $A$ is a $\sigma$-PBW extension of a ring $R$, and $M_R$ is a skew-Armendariz module with $R\subseteq M$, then $M_R$ is an abelian module.
\begin{proof}
The assertion follows from Theorem \ref{generAlhevaz2012Theorem2.26} and the fact that every skew-Armendariz module is linearly skew-Armendariz module.
\end{proof}
\end{corollary}
Note that \cite{AlhevazMoussavi2012}, Corollary 2.27 is a particular case of Corollary \ref{generAlhevaz2012Corollary2.27}. Next theorem  generalizes \cite{AlhevazMoussavi2012}, Theorem 2.28.
\begin{theorem}\label{generAlhevaz2012Theorem2.28}
If $M_R$ is a reduced module, then $M_R$ is a p.p.-module if and only if $M_R$ is a p.q.-Baer module.
\begin{proof}
Consider $M_R$ a reduced right $R$-module. From Proposition \ref{generAlhevazMoussavi2012Lemma2.14}, we know that for every elements $m$ of $M$ and $r$ of $R$, the equality $mr=0$ implies $mRr=0$, which means that ${\rm ann}_R(\{m\}) = {\rm ann}_R(mR)$, and so, ${\rm ann}_R(\{m\}) = {\rm ann}_R(mR)$.
\end{proof}
\end{theorem}
Theorem \ref{generAlhevaz2012Theorem 2.29} generalizes  \cite{AlhevazMoussavi2012}, Theorem 2.29.
\begin{theorem}\label{generAlhevaz2012Theorem 2.29}
If $A$ is a $\sigma$-PBW extension of a ring $R$, and $M_R$ is an $(\Sigma, \Delta)$-compatible skew-Armendariz right $R$-module with $R\subseteq M$, then $M_R$ is p.p.-module if and only if $M\langle X\rangle_A$ is p.p.-module
\begin{proof}
Let $M_R$ be a p.p.-module. Consider an element $m$ of $M\langle X\rangle$ given by the expression $m=m_0 + m_1X_1 + \dotsb + m_kX_k$. We know that ${\rm ann}_R(\{m_i\}) = e_iR$, for idempotent elements $e_i\in R$, for every $i$. Consider the product of the elements $e$'s, that is, let $e:=e_0e_1\dotsb e_k$. Note that $e$ is idempotent ($M_R$ is abelian by Corollary \ref{generAlhevaz2012Corollary2.27}). Therefore we have $eR=\bigcap_{i=0}^{k} {\rm ann}_R(\{m_i\})$. From Theorem \ref{generAlhevaz2012Theorem2.25} we know that $\sigma_i(e)=e$ and $\delta_i(e)=0$, for every $i$, which guarantees that the product $me$ is zero, and so, $eA\subseteq {\rm ann}_A(\{m\})$. Now, if $g\in {\rm ann}_A(\{m\})$ is given by $g=b_0 + b_1X_1 + \dotsb + b_tX_t$, then $m_0b_j=0$ for $0\le j\le t$ ($M_R$ is skew-Armendariz). In this way, $b_0\in eR$, whence $g\in eA$, that is, ${\rm ann}_A(\{m\})=eA$. In other words, we have shown that $M\langle X\rangle$ is a p.p.-module over $A$.

Conversely, if $M\langle X\rangle_A$ is a p.p.-module and $m$ is an element of $M$, for an idempotent element $e\in A$ given by $e=e_0 + e_1X_1 + \dotsb + e_pX_p$, we have $e(1-e)=0=(1-e)e$, or equivalently, $(e_0 + e_1X_1 + \dotsb + e_pX_p)(1-e_0 - e_1X_1 - \dotsb - e_pX_p) = 0 = (1-e_0 - e_1X_1 - \dotsb - e_pX_p)(e_0 + e_1X_1 + \dotsb + e_pX_p)$. As we know, $M_R$ is skew-Armendariz, so $e_0(1-e_0)=0$ and $(1-e_0)e_i=0$, for every $i$, which means that $e_0e_i=0,\ e_i=e_0e_i$, that is, $e_i=0$. Then $e = e_0^{2}=e_0\in R$, and ${\rm ann}_A(\{m\}) = eA$, whence ${\rm ann}_R(\{m\}) = eR$, i.e., $M_R$ is a p.p.-module.
\end{proof}
\end{theorem}
Theorem \ref{generAlhevaz2012Theorem2.30} extends \cite{AlhevazMoussavi2012}, Theorem 2.30, from Ore extensions to $\sigma$-PBW extensions.
\begin{theorem}\label{generAlhevaz2012Theorem2.30}
If $A$ is a $\sigma$-PBW extension of a ring $R$, $M_R$ is an $(\Sigma, \Delta)$-compatible skew-Armendariz module with $R\subseteq M$, then $M_R$ is Baer if and only if $M\langle X\rangle_A$ is Baer.
\begin{proof}
Suppose that $M_R$ is a Baer module. Let $J\subseteq M\langle X\rangle$ and $J_0$ the set of elements $m$ of $M$ such that $m$ is the leading coefficient of some non-zero element of $J$. Using that $M_R$ is Baer, there exists $e^{2}=e\in R$ with ${\rm ann}_R(J_0)=eR$, and hence $eA\subseteq {\rm ann}_A(J)$, by Proposition \ref{generAlehevaz2012Lemma2.15}. Now, consider an element $g=b_0 + b_1X_1 + \dotsb + b_tx_t\in {\rm ann}_A(J)$. Since $M_R$ is skew-Armendariz, $J_0b_j=0$, for $0\le j\le t$. This fact means that $b_j=eb_j$, for every $j$, and $g=eg\in A$, so ${\rm ann}_A(J) = EA$ and $M\langle X\rangle_A$ is a Baer module. Finally, if $M\langle X\rangle_A$ is a Baer module y $C\subseteq M$, then $C[X]\subseteq M\langle X\rangle$, and since $M\langle X\rangle$ is Baer, there exists an idempotent element $e = e_0 + e_1X_1 + \dotsb + e_pX_p\in A$ with ${\rm ann}_A)[C[X]] = eA$. In this way, $Ce_0=\{0\}$ and $e_0R\subseteq {\rm ann}_R(C)$. On the other hand, if $r\in {\rm ann}_R(C)$, then $C[X]r=0$ (Proposition \ref{generAlhevaz2012Lemma2.16}), whence $t=et$, that is, $t=e_0t\in e_0R$, which proves that ${\rm ann}_R(C)=e_0R$, and hence $M_R$ is a Baer module.
\end{proof}
\end{theorem}
The next proposition generalizes \cite{AlhevazMoussavi2012}, Proposition 2.32, Corollaries 2.33 and 2.34.
\begin{proposition}\label{generAlhevaz2012Proposition2.32}
Let $A$ be a $\sigma$-PBW extension of a ring $R$, and let $M_R$ be an $(\Sigma, \Delta)$-compatible and reduced module. If $m$ is a torsion element in $M\langle X\rangle$, i.e., $mh=0$, for some non-zero element $h\in A$, then there exists a non-zero element $c\in R$ such that $mc=0$.
\begin{proof}
Consider the elements $m=m_0 + m_1X_1 + \dotsb + m_kX_k\in M\langle X\rangle,\ 0\neq f=b_0 + b_1Y_1 + \dotsb + a_tY_t\in A,\ a_t\neq 0$, with $mf=0$. We have $mf= (m_0+m_1X_1+\dotsb + m_kX_k)(b_0+b_1Y_1+\dotsb + b_tY_t)=\sum_{l=0}^{k+t} \biggl(\sum_{i+j=l} m_iX_ib_jY_j\biggr)$. Note that ${\rm lc}(mf)= m_k\sigma^{\alpha_k}(b_t)c_{\alpha_k, \beta_t}=0$. Since $A$ is bijective, $m_k\sigma^{\alpha_k}(b_t)=0$, and by the $\Sigma$-compatibility of $M_R$, $m_kb_t=0$. The idea is to prove that $m_pb_q=0$ for $p+q\ge 0$.  We proceed  by induction. Suppose that $m_pb_q=0$ for $p+q=k+t, k+t-1, k+t-2, \dotsc, l+1$ for some $l>0$. By Proposition \ref{generAlhevazMoussavi2012Lemma2.14} and expression (\ref{perooo}), we obtain $m_pX_pb_qY_q=0$, for these values of $p+q$. In this way, we only consider the sum of the products $m_uX_ub_vY_v$, where $u+v=l, l-1,l-2,\dotsc, 0$. Fix $u$ and $v$. Consider the sum of all terms of $mf$  having exponent $\alpha_u+\beta_v$. From expression (\ref{perooo}), Proposition \ref{generAlhevazMoussavi2012Lemma2.14}, and the assumption $mf=0$, we know that the sum of all coefficients of all these terms  can be written as
\begin{equation}\label{Federer}
m_u\sigma^{\alpha_u}(b_v)c_{\alpha_u, \beta_v} + \sum_{\alpha_{u'} + \beta_{v'} = \alpha_u + \beta_v} m_{u'}\sigma^{\alpha_{u'}} ({\rm \sigma's\ and\ \delta's\ evaluated\ in}\ b_{v'})c_{\alpha_{u'}, \beta_{v'}} = 0.
\end{equation}
By assumption, we know that $m_pb_q=0$, for $p+q=k+t, k+t-1, \dotsc, l+1$.  So, Proposition \ref{generAlhevazMoussavi2012Lemma2.14} guarantees that the product
\[m_p({\rm \sigma's\ and\ \delta's\ evaluated\ in}\ b_{q})\ \ \ \ \ \ \ ({\rm any\  order\ of}\ \sigma's\ {\rm and}\ \delta's)
\]
is equal to zero. Then  $[({\rm \sigma's\ and\ \delta's\ evaluated\ in}\ b_{q})a_p]^2=0$, and hence we obtain the equality $({\rm \sigma's\ and\ \delta's\ evaluated\ in}\ b_{q})m_p=0$ ($M_R$ is reduced). In this way, multiplying (\ref{Federer}) by $m_l$, and using the fact that the elements $c_{i,j}$ in Definition \ref{gpbwextension} (iv) are in the center of $R$,
\begin{equation}\label{doooooo}
m_u\sigma^{\alpha_u}(b_v)a_kc_{\alpha_u, \beta_v} + \sum_{\alpha_{u'} + \beta_{v'} = \alpha_u + \beta_v} m_{u'}\sigma^{\alpha_{u'}} ({\rm \sigma's\ and\ \delta's\ evaluated\ in}\ b_{v'})m_lc_{\alpha_{u'}, \beta_{v'}} = 0,
\end{equation}
whence, $m_u\sigma^{\alpha_u}(b_0)a_l=0$. Since $u+v=l$ and $v=0$, then $u=l$, so $m_l\sigma^{\alpha_l}(b_0)m_l=0$, i.e., $[m_l\sigma^{\alpha_l}(b_0)]^{2}=0$, from which $m_l\sigma^{\alpha_l}(b_0)=0$ and $m_lb_0=0$, by Proposition \ref{generAlhevazMoussavi2012Lemma2.14}.  Therefore, we now have to study the expression (\ref{Federer}) for $0\le u \le l-1$ and $u+v=l$. If we multiply (\ref{doooooo}) by $m_{l-1}$, we obtain
{\small{\[
m_u\sigma^{\alpha_u}(b_v)m_{l-1}c_{\alpha_u, \beta_v} + \sum_{\alpha_{u'} + \beta_{v'} = \alpha_u + \beta_v} m_{u'}\sigma^{\alpha_{u'}} ({\rm \sigma's\ and\ \delta's\ evaluated\ in}\ b_{v'})m_{k-1}c_{\alpha_{u'}, \beta_{v'}} = 0.
\]}}
Using a similar reasoning as above, we can see that $m_u\sigma^{\alpha_u}(b_1)m_{l-1}c_{\alpha_u, \beta_1}=0$. Since $A$ is bijective, $m_u\sigma^{\alpha_u}(b_1)m_{l-1}=0$, and using the fact $u=l-1$, we have $[m_{l-1}\sigma^{\alpha_{l-1}}(b_1)]=0$, which imply $m_{l-1}\sigma^{\alpha_{l-1}}(b_1)=0$, that is, $m_{l-1}b_1=0$. Continuing in this way we prove that $m_ib_j=0$, for $i+j=l$. Therefore $a_ib_j=0$, for $0 \le i\le m$ and $0\le j\le t$. Finally, since $f$ is a non-zero element of $A$, we may consider $c:=a_t\neq 0$, and hence $mc=0$, by Proposition \ref{generAlhevaz2012Lemma2.16}.
\end{proof}
\end{proposition}
\subsection{Skew quasi-Armendariz modules over $\sigma$-PBW extensions}\label{SkewquasiArmendarizmodules}
In \cite{Hirano2002}, Hirano called a ring $B$ {\em quasi-Armendariz}, if whenever $f(x)B[x]g(x)=0$, where $f(x)=a_0 + a_1x + \dotsb + a_mx^{m}$ and $g(x)=b_0 + b_1x+\dotsb + b_nx^{n}\in B[x]$, then $a_iBb_j=0$, for all $i$ and $j$. In the same paper, Hirano called a right $B$-module {\em quasi-Armendariz}, if whenever $m(x)B[x]f(x)=0$, where $m(x)=m_0 + m_1x + \dotsb + m_sx^{s}\in M\langle X\rangle$ and $f(x)=a_0 + a_1x + \dotsb + a_tx^{t}\in B[x]$, implies that $m_iBa_j=0$, for all $i, j$. Now, in \cite{AlhevazMoussavi2012}, Definition 3.1, Alhevaz and Moussavi introduced the notion of {\em skew quasi-Armendariz} for modules in the following way: let $M_B$ be a right $B$-module, $\alpha$ an endomorphism of $B$ and $\delta$ an $\alpha$-derivation of $B$. $M_B$ is called {\em skew quasi-Armendariz}, if whenever $m(x)=m_0 + m_1x+\dotsb + m_kx^{k}\in M\langle X\rangle$, $f(x)=b_0+b_1x+\dotsb + b_nx^{n}\in B[x;\alpha, \delta]$ satisfy $m(x)B[x;\alpha, \delta]f(x)=0$, we have $m_ix^{i}Bx^{t}b_jx^{j} = 0$ for $t\ge 0,\ i=0,1,\dotsc, k$ and $j=0, 1,\dotsc, n$. With the aim of extending these definitions for the class of $\sigma$-PBW extensions, we present the following definition:
\begin{definition}
If $A$ is a $\sigma$-PBW extension of a ring $R$, and $M_R$ is a right $R$-module, $M_R$ is called {\em skew quasi-Armendariz}, if whenever $m=m_0 + m_1X_1+ \dotsb + m_kX_k\in M\langle X\rangle,\ f=a_0 + a_1Y_1 + \dotsb + a_pY_p\in A$ with $mAf=0$, then $m_iX_iRX_ta_jY_j = 0$, for any $X_t\in {\rm Mon}(A)$, and any values of $i, j$.
\end{definition}
The next theorem generalizes \cite{AlhevazMoussavi2012}, Theorem 3.2, for the Ore extensions to $\sigma$-PBW extensions, and its proof is easy.
\begin{theorem}
If $A$ is a quasi-commutative $\sigma$-PBW extension of $R$, and $M_R$ is an $\Sigma$-compatible right $R$-module, then:
\begin{enumerate}
\item [\rm (1)] The following assertions are equivalent:
\begin{enumerate}
\item [\rm (a)] for every $m\in M\langle X\rangle_A$, $({\rm ann}_A\{mA\}\cap R)\langle x_1,\dotsc, x_n\rangle = {\rm ann}_A\{mA\}$.
\item [\rm (b)] for every $m=m_0 + m_1X_1 + \dotsb + m_kX_k\in M\langle X\rangle_A$ and $f=a_0 + a_1Y_1 + \dotsb + a_tY_t\in A$, the equality $mAf=0$ implies $m_iRa_j=0$, for every $i, j$.
\end{enumerate}
\item [\rm (2)] If $M_R$ is a skew quasi-Armendariz module and $m$ is an element of $M\langle X\rangle_A$, the equality ${\rm ann}_A\{mA\}\neq \{0\}$ implies that ${\rm ann}_A\{mA\}\cap R\neq \{0\}$.
\end{enumerate}
\end{theorem}
Theorem \ref{generAlhevazMoussavi2012} generalizes \cite{AlhevazMoussavi2012}, Theorem 3.3.
\begin{theorem}\label{generAlhevazMoussavi2012}
If $A$ is a quasicommutative $\sigma$-PBW extension of $R$, and $M_R$ is an $\Sigma$-compatible right $R$-module, then the following assertions are equivalent:
\begin{enumerate}
\item [\rm (1)] $M_R$ is a skew quasi-Armendariz module;
\item [\rm (2)] The map $\psi:{\rm Ann}_R({\rm sub}(M_R)) \to {\rm Ann}_A({\rm sub}(M\langle X\rangle_A))$, defined by $\psi({\rm ann}_R(M)) = {\rm ann}_A(N)$ \linebreak $={\rm ann}_A(N[X])$, for all $N\in {\rm sub}(M_R)$, is bijective, where ${\rm sub}(M_R)$ and ${\rm sub}(M\langle X\rangle_A)$ denote the sets of submodules of $M_R$ and $M\langle X\rangle_A$, respectively.
\end{enumerate}
\begin{proof}
The proof follows from \cite{ReyesSuarezUMA2017}, Theorem 4.12.
\end{proof}
\end{theorem}
Theorem \ref{generAlhevaz2012Theorem3.6} generalizes  \cite{AlhevazMoussavi2012}, Theorem 3.6, Corollary 3.7 and \cite{Hirano2002}, Corollary 3.3.
\begin{theorem}\label{generAlhevaz2012Theorem3.6}
If $A$ is a $\sigma$-PBW extension of $R$, and $M_R$ is a skew quasi-Armendariz right $R$-module which is $(\Sigma, \Delta)$-compatible, then $M_R$ satisfies the ascending chain condition on annihilators of submodules if and only if so does $M\langle X\rangle_A$.
\begin{proof}
Suppose that $M_R$ satisfies the ascending chain condition on annihilators of submodules. Consider the chain of annihilator of submodules of $M\langle X\rangle_A$ given by $I_1\subseteq I_2\subseteq \dotsb $. Then, there exist submodules $K_i$ of $M\langle X\rangle_A$ with ${\rm ann}_A(K_i) = I_i$, for $i=1, 2,\dotsc$, and satisfying the relations $K_1\supseteq K_2\supseteq \dotsb$. Consider the sets $M_i$ consisting of all elements of $K_i$, for every $i$. By assumption, $M$ is skew-Armendariz, so that $M_i$ is a submodule of $M$, for every $i$. Note that $M_i\supseteq M_{i+1}$, whence ${\rm ann}_R(M_1)\subseteq {\rm ann}_R(M_2)\subseteq \dotsb$. Since we suppose that $M_R$ satisfies the ascending chain condition on annihilators of submodules, then there exists $p\ge 1$ such that ${\rm ann}_R(M_i)\subseteq {\rm ann}_R(M_p)$, for every value $i\ge p$. Our objective is to show that ${\rm ann}_A(K_i) = {\rm ann}_A(K_p)$, for these values $i\ge p$. With this in mind, consider an element $f\in {\rm ann}_A(K_i)$ given by the expression $f=a_0 + a_1X_1 + \dotsb + a_mX_m$. Since $M$ is skew quasi-Armendariz, we have $M_ia_j=0$, for $0\le j\le m$. This fact implies that $M_pa_j=0$, for $0\le j\le m$, and hence Proposition \ref{generAlhevaz2012Lemma2.16} guarantees that $K_pf=0$. Therefore we conclude that ${\rm ann}_A(K_i) = {\rm ann}_A(K_p)$, for $i\ge p$, i.e., $M\langle X\rangle_A$ satisfies the ascending chain condition on annihilator of submodules.

Conversely, suppose that $M\langle X\rangle_A$ satisfies the ascending chain condition on annihilator of submodules, and consider the chain of annihilator of submodules of $M_R$ given by $J_1\subseteq J_2\subseteq \dotsb$. Then, there exist submodules $M_i$ of $M$ with ${\rm ann}_R(M_i) = J_i$, for $i=1, 2,\dotsc$, and $M_1\supseteq M_2\supseteq \dotsb$. Note that $M_i[X]$ is a submodule of $M\langle X\rangle$, $M_i[X]\supseteq M_{i+1}[X]$, and, ${\rm ann}_A(M_i[X])\subseteq {\rm ann}_A(M_{i+1}[X])$, for every value of $i$. Since $M\langle X\rangle_A$ satisfies the ascending chain condition on annihilator of submodules, there exists $q\ge 1$ with ${\rm ann}_A(M_i[X]) = {\rm ann}_A(M_q[X])$, for $i\ge q$. Using that $M$ is skew-Armendariz, it follows that ${\rm ann}_R(M_i) = {\rm ann}_R(M_q)$, for $i\ge q$.
\end{proof}
\end{theorem}
The next theorem generalizes \cite{AlhevazMoussavi2012}, Theorem 3.9, \cite{HashemiMoussavi2005}, Corollary 2.8, \cite{BirkenmeierKimPark2001}, Corollary 2.8, and \cite{HongKimKwak2000}, Theorems 12 and 15.
\begin{theorem}\label{generalTheorem3.9}
If $A$ is a $\sigma$-PBW extension of $R$ and $M_R$ is an $(\Sigma, \Delta)$-compatible right $R$-module, then $M_R$ is quasi-Baer {\rm(}respectively, p.q.-Baer{\rm)} if and only if $M\langle X\rangle_{A}$ is quasi-Baer {\rm(}respectively, p.q.-Baer{\rm)}. In this case, $M_R$ is skew quasi-Armendariz.
\begin{proof}
Suppose that $M_R$ is a quasi-Baer right $R$-module. Let us see that $M_R$ is skew quasi-Armendariz. With this aim, consider the product
\[
(m_0 + m_1X_1+\dotsb + m_kX_k)A(a_0 + a_1Y_1+\dotsb + a_pY_p) = 0,
\]
where $m_0 + m_1X_1+\dotsb + m_kX_k\in M\langle X\rangle$ and $a_0 + a_1Y_1+\dotsb + a_pY_p\in A$. In particular, if we only take coefficients in $R$, we have the expression $(m_0 + m_1X_1+\dotsb + m_kX_k)R(a_0 + a_1Y_1+\dotsb + a_pY_p) = 0$. Since for any $r\in R$ we have the expressions
\begin{align}
0 = &\ (m_0 + m_1X_1+\dotsb + m_kX_k)r(a_0 + a_1Y_1+\dotsb + a_pY_p)\label{madrid} \\
= &\ m_0ra_0 + m_0ra_1Y_1 + \dotsb + m_0ra_pY_p + m_1X_1ra_0 + m_1X_1ra_1Y_1 + \dotsb + m_1X_1ra_pY_p \notag \\
+ &\ \dotsb + m_kX_kra_0 + m_kX_kra_1Y_1 + \dotsb + m_kX_kra_pY_p \notag \\
= &\ m_0ra_0 + m_0ra_1Y_1 + \dotsb + m_0ra_pY_p + m_1\sigma^{\alpha_1}(ra_0)X_1 + m_1p_{\alpha_1, ra_0} + m_1\sigma^{\alpha_1}(ra_1)X_1Y_1  \notag \\
+ &\ m_1p_{\alpha_1, ra_1}Y_1 + \dotsb + m_1\sigma^{\alpha_1}(ra_p)X_1Y_p + \dotsb + m_1p_{\alpha_1, ra_p}Y_p + \dotsb + m_k\sigma^{\alpha_k}(ra_0)X_k \notag  \\
+ &\ m_kp_{\alpha_k, ra_0} + m_k\sigma^{\alpha_k}(ra_1)X_kY_1 + m_kp_{\alpha_k, ra_1}Y_1 + \dotsb + m_k\sigma^{\alpha_k}(ra_p)X_kY_p + m_kp_{\alpha_k, ra_p}Y_p \notag \\
= &\ m_0ra_0 + m_0ra_1Y_1 + \dotsb + m_0ra_pY_p + m_1\sigma^{\alpha_1}(ra_0)X_1 + m_1p_{\alpha_1, ra_0}  \notag \\
+ &\ m_1\sigma^{\alpha_1}(ra_1)c_{\alpha_1,\beta_1}x^{\alpha_1+\beta_1} + m_1\sigma^{\alpha_1}(ra_1)p_{\alpha_1,\beta_1} + m_1p_{\alpha_1, ra_1}Y_1  \notag \\
+ &\ \dotsb + m_1\sigma^{\alpha_1}(ra_p)c_{\alpha_1,\beta_p}x^{\alpha_1+\beta_p} +  m_1\sigma^{\alpha_1}(ra_p)p_{\alpha_1,\beta_p} + \dotsb + m_1p_{\alpha_1, ra_p}Y_p  \notag \\
+ &\ \dotsb + m_k\sigma^{\alpha_k}(ra_0)X_k
+ m_kp_{\alpha_k, ra_0} + m_k\sigma^{\alpha_k}(ra_1)c_{\alpha_k, \beta_1}x^{\alpha_k+\beta_1} + m_k\sigma^{\alpha_k}(ra_1)p_{\alpha_k, \beta_1} \notag \\
+ &\ m_kp_{\alpha_k, ra_1}Y_1
+ \dotsb + m_k\sigma^{\alpha_k}(ra_p)c_{\alpha_k, \beta_p}x^{\alpha_k+\beta_p} + m_k\sigma^{\alpha_k}(ra_p)p_{\alpha_k,\beta_p} + m_kp_{\alpha_k, ra_p}Y_p, \notag
\end{align}
we can see that the leading coefficient of this product is $m_k\sigma^{\alpha_k}(ra_p)c_{\alpha_k, \beta_p}=0$, whence we obtain $m_k\sigma^{\alpha_k}(ra_p)=0$ (recall that the elements $c_{\alpha_k, \beta_p}$ are both invertible), and hence $m_kra_p=0$ because $M_R$ is $\Sigma$-compatible. In this way, $a_p\in {\rm ann}_R(m_kR)$. Now, since
\begin{align}
m_kX_krX_ta_pY_p = &\ (m_k\sigma^{\alpha_k}(r)X_k + m_kp_{\alpha_k, r}) (\sigma^{\alpha_t}(a_p)X_tY_p + p_{\alpha_t, a_p}Y_p) \notag \\
= &\ m_k\sigma^{\alpha_k}(r)X_k\sigma^{\alpha_t}(a_p)X_tX_p + m_k\sigma^{\alpha_k}(r)X_kp_{\alpha_t, a_p}Y_p \notag  \\
+ &\ m_kp_{\alpha_k,r}\sigma^{\alpha_t}(a_p)X_tX_p + m_kp_{\alpha_k,r}p_{\alpha_t, a_p}Y_p \notag \\
= &\ m_k\sigma^{\alpha_k}(r)(\sigma^{\alpha_k}(\sigma^{\alpha_t}(a_p))X_k + p_{\alpha_k, \sigma^{\alpha_t}(a_p)})X_tX_p + m_k\sigma^{\alpha_k}(r)X_kp_{\alpha_t, a_p}Y_p \notag  \\
+ &\ m_kp_{\alpha_k,r}\sigma^{\alpha_t}(a_p)X_tX_p + m_kp_{\alpha_k,r}p_{\alpha_t, a_p}Y_p \notag \\
= &\ m_k\sigma^{\alpha_k}(r\sigma^{\alpha_t}(a_p))X_kX_tX_p + m_k\sigma^{\alpha_k}(r)p_{\alpha_k, \sigma^{\alpha_t}(a_p)}X_tX_p  + m_k\sigma^{\alpha_k}(r)X_kp_{\alpha_t, a_p}Y_p \notag \\
+ &\ m_kp_{\alpha_k,r}\sigma^{\alpha_t}(a_p)X_tX_p +
m_kp_{\alpha_k,r}p_{\alpha_t, a_p}Y_p\label{force}
\end{align}
from Proposition \ref{generAlehevaz2012Lemma2.15} we obtain that all these expressions are zero, so $m_kX_krX_ta_pY_p=0$. As we know, $M_R$ is quasi-Baer, which means that there exists an element $e_k^{2}=e_k\in R$ with ${\rm ann}_R(m_kR)=e_kR$, whence $a_p=e_ka_p$. If we replace the element $r$ by $re_k$ in (\ref{madrid}), we obtain
\[
0 = (m_0 + m_1X_1+\dotsb + m_{k-1}X_{k-1})re_k(a_0 + a_1Y_1+\dotsb + a_pY_p),
\]
and using a similar reasoning to the above, we see that
\begin{align*}
0 = &\ (m_0 + m_1X_1+\dotsb + m_{k-1}X_{k-1})re_k(a_0 + a_1Y_1+\dotsb + a_pY_p)\\
= &\ m_0re_ka_0 + m_0re_ka_1Y_1 + \dotsb + m_0re_ka_pY_p + m_1\sigma^{\alpha_1}(re_ka_0)X_1 + m_1p_{\alpha_1, re_ka_0}  \notag \\
+ &\ m_1\sigma^{\alpha_1}(re_ka_1)c_{\alpha_1,\beta_1}x^{\alpha_1+\beta_1} + m_1\sigma^{\alpha_1}(re_ka_1)p_{\alpha_1,\beta_1} + m_1p_{\alpha_1, re_ka_1}Y_1  \notag \\
+ &\ \dotsb + m_1\sigma^{\alpha_1}(re_ka_p)c_{\alpha_1,\beta_p}x^{\alpha_1+\beta_p} +  m_1\sigma^{\alpha_1}(re_ka_p)p_{\alpha_1,\beta_p} + \dotsb + m_1p_{\alpha_1, re_ka_p}Y_p  \notag \\
+ &\ \dotsb + m_{k-1}\sigma^{\alpha_{k-1}}(re_ka_0)X_{k-1}
+ m_{k-1}p_{\alpha_{k-1}, re_ka_0} \\
+ &\ m_{k-1}\sigma^{\alpha_{k-1}}(re_ka_1)c_{\alpha_{k-1},
\beta_1}x^{\alpha_{k-1}+\beta_1}
+ m_{k-1}\sigma^{\alpha_{k-1}}(re_ka_1)p_{\alpha_{k-1}, \beta_1} \\
+ &\ m_{k-1}p_{\alpha_{k-1}, re_ka_1}Y_1
+ \dotsb + m_{k-1}\sigma^{\alpha_{k-1}}(re_ka_p)c_{\alpha_{k-1}, \beta_p}x^{\alpha_{k-1}+\beta_p}\\
+ &\ m_{k-1}\sigma^{\alpha_{k-1}}(re_ka_p)p_{\alpha_{k-1},\beta_p}
+ m_{k-1}p_{\alpha_{k-1}, re_ka_p}Y_p,
\end{align*}
whence $m_{k-1}\sigma^{\alpha_{k-1}}(re_ka_p)c_{\alpha_{k-1}, \beta_p}= 0\Rightarrow m_{k-1}\sigma^{\alpha_{k-1}}(re_ka_p) = 0\Rightarrow  m_{k-1}re_ka_p = 0$, which implies that $m_{k-1}Ra_p=0$, i.e., $a_p\in {\rm ann}_R(m_{k-1}R)$. Again, using a similar reasoning to the above in expression (\ref{force}), we can see that $m_{k-1}X_{k-1}RX_ta_pX_p=0$. Therefore we have showed that $a_p\in {\rm ann}_R(m_kR)\cap {\rm ann}_R(m_{k-1}R)$. Since $M_R$ is quasi-Baer, there exists an idempotent element $s\in R$ with ${\rm ann}_R(m_{k-1}R)=sR$, whence $a_p=sa_p$. If we define the element $e_{k-1}$ as $e_{k-1}=e_{k}s$, then we can see that $e_{k-1}\in {\rm ann}_R(m_kR)\cap {\rm ann}_R(m_{k-1}R)$. Now, if we replace $r$ by $re_{k-1}$ in (\ref{madrid}), we obtain the equality
\[
0 = (m_0 + m_1X_1+\dotsb + m_{k-2}X_{k-2})re_{k-1}(a_0 + a_1Y_1+\dotsb + a_pY_p).
\]
We can show that the relations $m_{k-2}\sigma^{\alpha_{k-2}}(re_{k-1}a_p)c_{\alpha_{k-2}, \beta_p}= 0\Rightarrow m_{k-2}\sigma^{\alpha_{k-2}}(re_{k-1}a_p) = 0\Rightarrow  m_{k-2}re_{k-1}a_p = 0$, which implies that $m_{k-2}Ra_p=0$, i.e., $a_p\in {\rm ann}_R(m_{k-2}R)$, and hence $m_{k-2}X_{k-2}RX_ta_pX_p=0$. Continuing in this way, we can prove that $m_iX_iRX_ta_pY_p=0$, for $i=0,\dotsc, k$, and any $X_t\in {\rm Mon}(A)$. Similarly, using the total order on ${\rm Mon}(A)$, where $Y_p\succ Y_{p-1}\succ \dotsb Y_1 \succ 1$, we can show that $m_iX_iRX_ta_{p-1}Y_{p-1} = m_iX_iRX_ta_{p-2}Y_{p-2} = \dotsb = m_iX_iRX_ta_1Y_1 = m_iX_iRX_ta_0=0$, which allows us to conclude that $M_R$ is skew quasi-Armendariz.

Now, we will prove that $M\langle X\rangle_A$ is quasi-Baer. Let $J$ be a $A$-submodule of $M\langle X\rangle$, and consider the set $N$ as the union of the set of the leading coefficients of non-zero elements of $J$ with the set $\{0\}$. Note that $N$ is a submodule of $M$. By assumption, $M_R$ is quasi-Baer, so there exists an idempotent element $e$ of $R$ with ${\rm ann}_R(N)=eR$, which implies that $eA\subseteq {\rm ann}_A(J)$ (Proposition \ref{generAlehevaz2012Lemma2.15}). With the aim of proving that $eA\supseteq {\rm ann}_A(J)$, consider an element $f=a_0 + a_1x_1 + \dotsb + a_pX_p\in {\rm ann}_A(J)$. Since $M_R$ is skew quasi-Armendariz, it follows that  $Na_j=0$, for $0\le j \le p$. Then, $b_j=eb_j$, for every $j$, an $f=ef\in eA$, which guarantees that $eA\supseteq {\rm ann}_A(J)$, and hence, $eA = {\rm ann}_A(J)$.
Conversely, if $M\langle X\rangle_A$ is quasi-Baer and $I$ is a submodule of $M$, it follows that $I[X]$ is a submodule of $M\langle X\rangle$, and since $M\langle X\rangle$ is quasi-Baer, there exists $e^{2}=e = e_0+e_1+\dotsb + e_lX_l\in A$ with ${\rm ann}_A(I[X]) = eA$. Note that $Ie_0=0$ and $e_0R\subseteq {\rm ann}_R(I)$. Finally, if $s\in {\rm ann}_R(I)$, then $I[X]s=0$ (Proposition \ref{generAlhevaz2012Lemma2.16}), and therefore $t=et$, which implies $t=e_0t\in e_0R$, that is, $e_0R\supseteq {\rm ann}_R(I)$, i.e., $e_0R= {\rm ann}_R(I)$. This concludes the proof.
\end{proof}
\end{theorem}
\begin{remark}\label{Notpassp.q.Baer}
\begin{itemize}
\item Note that a ring $B$ is right p.q.-Baer if and only if $B_B$ is a p.q.-Baer module. However, this does not hold for the property of being p.q.-Baer. More exactly, there exists a p.q.-Baer right $B$-module such that $B$ is not right p.q.-Baer (\cite{AlhevazMoussavi2012}, Example 3.10).
\item The condition on the $(\Sigma, \Delta)$-compatibility in Theorem \ref{generalTheorem3.9} can not be dropped, since there exists an example of a ring $B$ such that $B[x;\delta]$ is Baer, and hence quasi-Baer, but $B$ is not quasi-Baer, see \cite{Armendariz1974}, Example 1.
\end{itemize}
\end{remark}

\vspace{0.5cm}

\noindent {\bf \Large{Acknowledgements}}

\vspace{0.5cm}

The author was supported by the research fund of Departamento de Matem\'aticas, Universidad Nacional de Colombia, Bogot\'a, Colombia, HERMES CODE 30366.


\end{document}